\documentclass[12pt]{article}
\usepackage{amssymb,amsmath,amsopn,amsfonts}
\usepackage[dvips]{color}
\usepackage{colordvi,multicol}
\usepackage{epsf}
\newfam\msbfam
\font\tenmsb=msbm10 \textfont\msbfam=\tenmsb \font\sevenmsb=msbm7
\scriptfont\msbfam=\sevenmsb \font\fivemsb=msbm5
\scriptscriptfont\msbfam=\fivemsb

\def\th#1{\vspace{1mm}\noindent{\bf #1}\quad}

\def\proof{\vspace{1mm}\noindent{\it Proof}\quad}

\numberwithin{equation}{section}

\def\bc{\begin{center}}
\def\ec{\end{center}}
\def\no{\noindent}
\def\hang{\hangindent\parindent}
\def\textindent#1{\indent\llap{\qquad #1\ \ \enspace}\ignorespaces}
\def\ref{\par\hang\textindent}

\begin{document}
\centerline{\Large\bf BSDE and generalized Dirichlet forms:}
\centerline{\Large\bf the infinite
dimensional case}

\vspace{0.5 true cm}

 \centerline{RONG-CHAN ZHU$^{\textrm{a,b}}$} \centerline{\small
a. Institute of Applied Mathematics, Academy of Mathematics and Systems
Science,} \centerline{\small Chinese Academy of Sciences, Beijing
100190, China}
\centerline{\small b. Department of Mathematics, University of Bielefeld, D-33615 Bielefeld, Germany}
\centerline{\small E-mail: zhurongchan@126.com}

\footnotetext{\footnotesize Research supported by 973 project, NSFC, key Lab of CAS, the DFG through IRTG 1132 and CRC 701}
 \vskip
1.4cm
\begin{abstract}
\vskip 0.1cm \noindent We consider the following
quasi-linear parabolic system of backward partial differential equations on a Banach space $E$
 $$(\partial_t+L)u+f(\cdot,\cdot,u, A^{1/2}\nabla u)=0 \textrm{ on } [0,T]\times E,\qquad u_T=\phi,$$
 where $L$ is a possibly degenerate  second order differential operator with merely measurable coefficients. We solve this system in the framework of generalized Dirichlet forms and employ the stochastic calculus
 associated to the Markov process with generator $L$ to obtain a probabilistic representation of the solution $u$ by solving the corresponding backward stochastic differential equation.
The solution satisfies the corresponding mild equation which is equivalent to being a generalized solution of the PDE. A further main result is the generalization of  the martingale representation theorem in infinite dimension
using  the stochastic calculus associated to the generalized Dirichlet form given by $L$. The nonlinear term $f$ satisfies a monotonicity condition with respect to $u$ and a Lipschitz condition with respect to $\nabla u$.
\end{abstract}
\no{\footnotesize{\bf 2010 Mathematics Subject Classification}: 31C25, 35R15, 60H30}

\no{\footnotesize{\bf Keywords}:\hspace{2mm}  backward stochastic differential equations, quasi-linear parabolic partial differential equations, Dirichlet
forms, generalized Dirichlet forms, Markov processes, martingale representation, Kolmogorov equations}

\section{Introduction}
Consider the following quasi-linear parabolic system of backward partial differential equations on a (real)  Banach space $E$
 $$(\partial_t+L)u+f(\cdot,\cdot,u, A^{1/2}\nabla u)=0 \textrm{ on } [0,T]\times E,\qquad u_T=\phi, \eqno{(1.1)}$$
 where $L$ is a second order differential operator with measurable coefficients, $\nabla u$ is the $H$-gradient of $u$ and $(H,\langle \cdot, \cdot \rangle_H)$ is
 a separable real Hilbert space such that $H\subset E$ densely and continuously. $A$ is a symmetric, positive-definite and bounded operator on $H$. This equation is also called nonlinear  Kolmogorov equation on an infinite dimensional space. In fact, in this paper we study systems of PDE
of type (1.1), i.e. $u$ takes values in $\mathbb{R}^l$ for some fixed $l\in \mathbb{N}$. For simplicity, in this introductory section we explain our results in the case $l=1$.

 Various concepts of solution are known for (linear and) nonlinear parabolic equations in infinite dimensions. In this paper we will consider solutions in the  sense of Definition 2.4, i.e.
there is a sequence  $\{u^n\}$ of strong solutions with data $(\phi^n,f^n) $ such that $$\|u^n-u\|_T\rightarrow0,\|\phi^n-\phi\|_2\rightarrow0 \textrm{ and } \lim_{n\rightarrow\infty} f^n=f \textrm{ in } L^1([0,T];L^2).$$
We will prove the above definition for  solution is equivalent to being a solution of the following mild equation in $L^2$ sense
$$u(t,x)=P_{T-t}\phi(x)+\int_t^TP_{s-t}f(s,\cdot,u_s,A^{1/2} \nabla u_s)(x)ds.\eqno{(1.2)}$$ This formula is meaningful provided $u$ is even only once differentiable with respect to $x$. Thus, the solutions we consider are in a sense intermediate between classical and viscosity solutions.

The notion of viscosity solution, developed by many authors, in particular M. Crandall and P. L. Lions and their collaborators, is not discussed here. Generally speaking, the class of equations that can be treated by this method (c.f. [22-24] ) is much more general than those considered in this paper: it includes fully nonlinear operators. However, none of these results  are  applicable to our situation because we only need that the coefficients of  the operator $L$ are measurable.

If $E$ is a Hilbert space, in \cite{FG}, mild solutions  of the above PDE have been considered, and  a probabilistic technique, based on backward stochastic differential equations, has been used
 to prove the existence and uniqueness for  mild solutions. Furthermore, their results have been extended in \cite{BC} and \cite{M}. All these results  need some regular conditions  for the coefficients of $L$ and $f$ to make sure that the process $X$ has regular dependence on parameters, which are not needed for  our results. In this  paper, we
 will prove existence and uniqueness of a  solution $u$ by methods  from functional analysis. In fact this paper is an extension of our paper \cite{Z} to the infinite dimensional case. Though \cite{Z} serves as a guideline, serious obstacles appear at various places if
$E$ is infinite dimensional, which we overcome in this work (see e.g. Proof of Theorem 3.8).

The connection between backward stochastic equations and nonlinear partial differential equations was proved for the finite dimensional case
 e.g. in  \cite{BPS05},\cite{BDHPS}, \cite{PP}  ( see also the  references therein).  A further  motivation of this paper  is to give a probabilistic interpretation for  the solutions of the above PDE's, i.e. in this infinite dimensional case.

If $E$ is equal to a Hilbert space $H$, and the coefficients of the second-order differential operator $L$ are sufficiently regular,  then PDE (1.1)
  has a classical solution and one may construct the pair of processes $Y_s^{t,x}:=u(s,X_s^{t,x}), Z_s^{t,x}:=A^{1/2}\nabla u(s,X_s^{t,x})$ where $X_s^{t,x}, t\leq s\leq T$,
 is the diffusion process with infinitesimal operator $L$ which starts from $x$ at time  $t$. Then, using It\^{o}'s formula one checks that $(Y_s^{t,x}, Z_s^{t,x})_{t\leq s\leq T}$ solves the BSDE
$$Y_s^{t,x}=\phi(X_T^{t,x})+\int_s^Tf(r,X_r^{t,x},Y_r^{t,x},Z_r^{t,x})dr-\int_s^T\langle Z_r^{t,x},dW_r\rangle_H,\eqno{(1.3)}$$
Here $W_r$ is a cylindrical Wiener process in $H$. Conversely, for regular coefficients  by standard methods one can prove that (1.3) has a unique solution $(Y_t^{s,x}, Z_t^{s,x})_{s\leq t\leq T}$ and then $u(s,x):=Y_s^{s,x}$ is a solution to PDE (1.1).
If $f$ and the coefficient of $L$ are  Lipschitz continuous then in  \cite{FG} the authors  prove the probabilistic interpretation  above remains true, if one considers mild solutions to PDE (1.1). There are many papers
that study  forward-backward systems in  infinite dimension (cf \cite{FG}, \cite{FH} and the references therein). In
 these approaches, since the coefficients are Lipschitz continuous, the Markov process $X$ with infinitesimal operator $L$ is a diffusion process which satisfies an SDE and so one can use its associated stochastic calculus.

In \cite{BPS05} Bally, Pardoux and Stoica consider a semi-elliptic symmetric second-order differential operator $L$
( which is written in divergence form ) with measurable coefficients in finite dimension. They prove the above system of PDE has a unique solution $u$ in some functional space.
 Then they prove the solution $Y^{t.x}$ of the BSDE yields a precised version of the  solution $u$ so that one has $Y_s^{t,x}=u(t+s,X_s), P^x$-a.s. In this paper, we generalize their results to  a non-symmetric second order differential operator $L$ in infinite dimensions.

In this paper, we consider  PDE (1.1)  for a non-symmetric second order differential operator $L$ in infinite dimensions, which is associated to the bilinear form
$$\mathcal{E}(u,v)=\int \langle A(z)\nabla u(z),\nabla v(z)\rangle_H d\mu(z)+\int \langle A(z)b(z),\nabla u(z)\rangle_H v(z) d\mu(z), u,v \in \mathcal{F}C_b^\infty.$$
Here we only need $|A^{1/2}b|_H \in L^2(E;\mu)$. That is to say, in general the above bilinear form $\mathcal{E}$ does not satisfy any  sector condition. We use the theory of generalized Dirichlet form and
 the associated stochastic calculus( cf.  [35-38]) to generalize the results in \cite{BPS05} both to infinite dimensional state spaces and fully nonsymmetric operators $L$.

In the analytic part of our paper, we don't need $\mathcal{E}$ to be a generalized Dirichlet form. We start from a semigroup $(P_t)$ satisfying conditions (A1)-(A3), specified in Section 2 below (see, in particular, also Remark 2.1 (viii)).
Such a semigroup can e.g.  be constructed from a generalized Dirichlet form. It can also be constructed by other methods (see e.g. \cite{DR}). Under conditions (A1)-(A3), the coefficients of $L$ may be quite singular and only very broad assumptions on $A$ and $b$ are needed.

The paper is organized as follows. In Sections 2 and 3, we use functional analytic methods to solve  PDE (1.1)  in the sense of Definition 2.4 or equivalently in the sense of (1.2).
Here the function $f$ need not be Lipschitz continuous with respect to $y$; monotonicity suffices.  And $\mu$, which appears in the monotonicity
 conditions (see condition (H2) in Section 3.2 below), can depend on $t$. $f$ is, however, assumed to be Lipschitz continuous with respect to the last variable.
We emphasize that the first order term of $L$ with coefficient $Ab$   cannot be incorporated into $f$ unless it  is bounded. Hence we are forced to take it as a part of $L$ and
have to consider a diffusion process $X$  which is generated by an operator $L$ which is the generator of a (in general non-sectorial) generalized Dirichlet form. We also emphasize that under our conditions PDE (1.1)
 cannot be tackled by standard monotonicity methods (see e.g. \cite{Ba10}) because of  lack of a suitable Gelfand triple $\mathcal{V}\subset \mathcal{H }\subset \mathcal{V}^*$ with $\mathcal{V}$ being a reflexive Banach space.

In Section 4,  we assume that $\mathcal{E}$ is a generalized Dirichlet form and is associated with a strong Markov process $X=(\Omega,\mathcal{F}_\infty,\mathcal{F}_t,X_t,P^x)$. Such a process can be constructed if $\mathcal{E}$ is quasi-regular. We extend the stochastic calculus for the Markov process  in order to generalize the martingale
 representation theorem.  More precisely, in order to treat BSDE's, in Theorem 4.3 we show   that there is  a set  $\mathcal{N}$
of null capacity outside of which the following representation theorem holds : for  every bounded $\mathcal{F}_\infty$-measurable random variable $\xi$, there  exists a predictable process $\phi:[0,\infty)\times \Omega\rightarrow H$, such that for each probability measure $\nu$, supported by $E\setminus\mathcal{N}$, one has
$$\xi=E^\nu(\xi|\mathcal{F}_0)+\sum_{i=0}^\infty\int_0^\infty\phi_s^idM_s^{i}\qquad P^\nu-a.e..$$
In fact, one may choose the exceptional set $\mathcal{N}$ such that if the process $X$ starts from a point of $\mathcal{N}^c$, it remains always in $\mathcal{N}^c$.
 As a consequence we deduce the existence of solutions for the BSDE using the existence of solutions  for PDE (1.1) in the usual way, however, only under $P^\mu$, because of our very  general coefficients of $L$ (c.f. Theorem 4.7).

 In Section 5, we employ the above results to deduce existence and uniqueness for the solutions of the BSDE under $P^x$ for $x\in \mathcal{N}^c$. As a consequence, in Theorem 5.4 one finds
 a version of the solution to PDE (1.1)  which satisfies the mild equation pointwise, i.e. for the solution $Y^{s}$ of the BSDE, we have $Y^{s}_t=u(t,X_{t-s}),P^x$-a.s. In particular, $Y_t^{t}$ is $P^x$-a.s.  equal to $u(t,x)$.

In Section 6, we give some examples of the operator $L$  satisfying our general conditions (A1)-(A5).

\section{The Linear Equation}

Let $E$ be a separable real Banach space and $(H,\langle \cdot, \cdot \rangle_H)$ a separable real Hilbert space such that $H\subset E$ densely and continuously.
Identifying $H$ with its topological dual $H'$ we obtain that $E'\subset H\subset E$ densely and continuously and $_{E'}\langle\cdot, \cdot\rangle_E=\langle\cdot,\cdot\rangle_H$ on $E'\times H$. Define the linear space of finitely based smooth functions on $E$ by
$$\mathcal{F}C_b^\infty:=\{f(l_1,...,l_m)| m\in N, f\in C_b^\infty(\mathbb{R}^m), l_1,..., l_m\in E'\}.$$
Here $C_b^\infty(\mathbb{R}^m)$ denotes the set of all infinitely differentiable (real-valued) functions with all partial derivatives bounded. For $u\in  \mathcal{F}C_b^\infty$ and $k\in E$ let
$$\frac{\partial u}{\partial k}(z):=\frac{d}{ds}u(z+sk)|_{s=0}, z\in E,$$
be the G\^{a}teaux derivative of $u$ in direction $k$. It follows that for $u=f(l_1,...,l_m)\in\mathcal{F}C_b^\infty$ and $k\in H$ we have that
$$\frac{\partial u}{\partial k}(z)=\sum_{i=1}^m\frac{\partial f}{\partial x_i}(l_1(z),...,l_m(z))\langle l_i, k\rangle_H,z\in E.$$
Consequently, $k\mapsto \frac{\partial u}{\partial k}(z)$ is continuous on $H$ and we can define $\nabla u(z)\in H$ by
$$\langle \nabla u(z), k\rangle_H=\frac{\partial u}{\partial k}(z).$$
Let $\mu$ be a finite positive measure on $(E,\mathcal{B}(E))$. By  $L_{sym}(H)$ we denote   the linear space of all symmetric and bounded linear operators on $H$
 equipped  with usual operator norm $\|\cdot\|_{L^\infty(H)}$. Let $A:E\mapsto L_{sym}(H)$ be measurable  such that $\langle A(z)h,h\rangle_H\geq0$ for all $z\in E, h\in H$ and let $b: E\rightarrow H$ be  $\mathcal{B}(E)/\mathcal{B}(H)$-measurable. Assume that the Pseudo inverse  $A^{-1}$ of $A$ is measurable

Denote the $H$-norm  by $|\cdot|_H$ and set  $\|u(z)\|_2^2:=\int|u(z)|^2d\mu(z)$ for $u \in L^2(E,\mu)$.
 We also denote $(u,v)_{L^2(E,\mu)}$ by $(u,v)$ for $u,v\in L^2(E,\mu)$. For $p\geq 1$, let $L^p(\mu)$, $L^p(\mu;H)$ denote $L^p(E,\mu)$, $L^p(E,\mu;H)$ respectively.  If $W$ is a function space, we will use $bW$ to denote
set of all  the bounded functions in $W$.

Furthermore,  we introduce the bilinear form
\begin{equation}\mathcal{E}(u,v):=\int \langle A(z)\nabla u(z),\nabla v(z)\rangle_H d\mu(z)+\int \langle A(z)b(z),\nabla u(z)\rangle_H v(z) d\mu(z), u,v \in \mathcal{F}C_b^\infty.\end{equation}
Consider the following conditions,

 \vskip.10in
\noindent (A1) $\langle A(\cdot)k,k\rangle\in L^1(\mu)$ and  bilinear form
 $$\mathcal{E}^A(u,v)=\int \langle A(z)\nabla u(z),\nabla v(z)\rangle_H d\mu(z); u,v \in \mathcal{F}C_b^\infty,$$ is closable on $L^2(E;\mu)$.
  \vskip.10in

 The closure of $\mathcal{F}C_b^\infty$ with respect to $\mathcal{E}^A_1:=\mathcal{E}^A+\langle\cdot,\cdot\rangle_H$ is denoted by $F$. Then $(\mathcal{E}^A,F)$
is a well-defined symmetric Dirichlet form on $L^2(E,\mu)$. We set $\mathcal{E}_1^A(u):=\mathcal{E}_1^A(u,u), u\in F.$
 \vskip.10in

\noindent(A2) Let $A^{1/2}b \in L^2(E;H,\mu)$, i.e. $\int|A^{1/2}b|_H^2d\mu<\infty$. We also assume there exists $\alpha\geq 0$ such that
\begin{equation}\int\langle Ab,\nabla u^2\rangle_H d\mu \geq -\alpha \|u\|_2^2, \qquad   u\in \mathcal{F}C_b^\infty, \forall u\geq 0.\end{equation}

\vskip.10in
Obviously, $\mathcal{E}$ from (2.1) immediately extends to all $u\in F, v\in bF$.
\vskip.10in

\noindent(A3) There exists a positivity preserving $C_0$-semigroup $P_t$ on $L^2(E;\mu)$ such that
 for any $t\in[0,T], \exists C_T>0$ such that $$\|P_t f\|_\infty \leq C_T \|f\|_\infty,$$ and such that  its $L^2$-generator  $(L,\mathcal{D}(L))$ has the following properties:  $b\mathcal{D}(L)\subset bF$ and
 for any $u\in bF$ there exists uniformly bounded $u_n\in \mathcal{D}(L)$ such that  ${\mathcal{E}}^{A}_1(u_n-u)\rightarrow0$ as $n\rightarrow\infty$  and that it is associated with the bilinear form $\mathcal{E}$ in (2.1)  in the sense that $\mathcal{E}(u,v)=-(Lu,v)$ for $u,v\in b\mathcal{D}(L)$.
 \vskip.10in
 To obtain a semigroup $P_t$ satisfying the above conditions, we can use generalized Dirichlet forms.
Let us recall the definition of a generalized Dirichlet form from \cite{St}.
 Let $E_1$ be a Hausdorff topological space and assume that its Borel
$\sigma$-algebra $\mathcal{B}(E_1)$ is generated by the set $C(E_1)$ of
all continuous functions on $E_1$. Let $m$ be a $\sigma$-finite
measure on $(E_1,\mathcal{B}(E_1))$ such that $\mathcal{H}:=L^2(E_1,m)$ is
a separable (real) Hilbert space. Let $(\mathcal{A},\mathcal{V})$ be
a coercive closed form on $\mathcal{H}$ in the sense of \cite{MR}. We  denote the
corresponding norm by $\|\cdot\|_{\mathcal{V}}$. Identifying
$\mathcal{H}$ with its dual $\mathcal{H}'$ we obtain that $
\mathcal{V}\rightarrow\mathcal{H}\cong\mathcal{H}'\rightarrow\mathcal{V}'$ densely and continuously.

Let $(\Lambda,D(\Lambda,\mathcal{H}))$ be a linear operator on
$\mathcal{H}$ satisfying the following assumptions:

(i) $(\Lambda,D(\Lambda,\mathcal{H}))$ generates a $C_0$-semigroup
of contractions $(U_t)_{t\geq0}$ on $\mathcal{H}$.

(ii)$\mathcal{V}$ is $\Lambda$-admissible, i.e.  $(U_t)_{t\geq}$
can be restricted to a $C_0$-semigroup on $\mathcal{V}$.

Let $(\Lambda,\mathcal{F})$ with corresponding norm $\|\cdot\|_{\mathcal{F}}$ be the closure of
$\Lambda:D(\Lambda,\mathcal{H})\cap\mathcal{V}\rightarrow\mathcal{V}'$ as an operator from $\mathcal{V}$ to $\mathcal{V}'$
 and $(\hat{\Lambda},\hat{\mathcal{F}})$ its dual operator.

 Let\begin{equation*}
\mathcal{E}(u,v)=\left\{\begin{array}{ll} \mathcal{A}(u,v)-\langle\Lambda u,v\rangle&\ \ \ \textrm{if} \ u\in \mathcal{F},v\in \mathcal{V},\\
\mathcal{A}(u,v)-\langle\hat{\Lambda} v,u\rangle&\ \ \ \ \textrm{if}
\ u\in \mathcal{V},v\in \hat{\mathcal{F}}.
\end{array}\right.
\end{equation*}
Here  $\langle\cdot,\cdot\rangle$ denotes the dualization between
$\mathcal{V}'$ and $\mathcal{V}$, which  coincides with the inner product $(\cdot,\cdot)_\mathcal{H}$ in $\mathcal{H}$ when restricted
to $\mathcal{H}\times \mathcal{V}$. We set
$\mathcal{E}_{\alpha}(u,v):=\mathcal{E}(u,v)+\alpha(u,v)_{\mathcal{H}}$
for $\alpha>0$. We call $\mathcal{E}$ the bilinear form associated
with $(\mathcal{A},\mathcal{V})$ and
$(\Lambda,D(\Lambda,\mathcal{H}))$. If $$u\in \mathcal{F}\Rightarrow
u^+\wedge1\in \mathcal{V} \textrm{ and }
\mathcal{E}(u,u-u^+\wedge1)\geq0,$$ then the bilinear form is called
a generalized Dirichlet form.  If the adjoint semigroup $(\hat{U}_t)_{t\geq0}$ of $(U_t)_{t\geq 0}$ can also be restricted to a $C_0$-semigroup on $\mathcal{V}$, let $(\hat{\Lambda},D(\hat{\Lambda},\mathcal{H}))$ denote the generator of $(\hat{U}_t)_{t\geq0}$ on $\mathcal{H}$, $\hat{\mathcal{A}}(u,v):=\mathcal{A}(v,u), u,v\in \mathcal{V}$ and let the coform $\hat{\mathcal{E}}$ be defined as the bilinear form associated with $(\hat{\mathcal{A}},\mathcal{V})$ and
$(\hat{\Lambda},D(\hat{\Lambda},\mathcal{H}))$.

 \vskip.10in

\th{Remark 2.1} (i) Some general criteria imposing conditions on $A$ and $\mu$ in order that $\mathcal{E}^A$ be closable
are e.g. given in  [25, Chap II, Section 2] and \cite{AR}.

(ii)In our case, due to our  general conditions on $b$ and $f$, we can't find  a suitable Gelfand triple $V\subset H\subset V^*$ with $V$ being a reflexive Banach space to apply  the monotonicity method as in \cite{Ba10} or  \cite{PR}.

(iii) We can construct a semigroup $P_t$ satisfying (A3) by the theory of generalized Dirichlet forms. More precisely, if  there exists a constant $\hat{c}\geq0$ such that $\mathcal{E}_{\hat{c}}(\cdot,\cdot):=\mathcal{E}(\cdot,\cdot)+\hat{c}(\cdot,\cdot)$ is a generalized Dirichlet form with domain $\mathcal{F}\times \mathcal{V}$ in one of the following three senses:

(a)$(E_1,\mathcal{B}(E_1),m)=(E,\mathcal{B}(E),\mu)$,

 $(\mathcal{A},\mathcal{V})=(\mathcal{E}^{A},F)$,

  $-\langle\Lambda u,v\rangle-\hat{c}(u,v)=\int \langle A(z)b(z),\nabla u(z)\rangle_H v(z) d\mu(z)$  for $u,v\in \mathcal{F}C_b^\infty$;

(b)$(E_1,\mathcal{B}(E_1),m)=(E,\mathcal{B}(E),\mu)$,

$\mathcal{A}\equiv0$ and $\mathcal{V}=L^2(E,\mu)$,

$-\langle\Lambda u,v\rangle=\mathcal{E}_{\hat{c}}(u,v)$ for $u,v\in D$, where $D\subset\mathcal{F}C_b^\infty$ densely w.r.t. $\mathcal{E}_1^A$-norm and $D\subset \mathcal{D}(L)$;

(c) $\mathcal{E}_{\hat{c}}=\mathcal{A}$, $\Lambda\equiv0$ (In this case $(\mathcal{E}_{\hat{c}},\mathcal{V})$ is a sectorial Dirichlet form in the sense of \cite{MR});

then there exists a sub-Markovian $C_0$-semigroup of contraction $P_t^{\hat{c}}$ associated with the generalized Dirichlet form $\mathcal{E}_{\hat{c}}$. Then $P_t:=e^{{\hat{c}}t}P_t^{\hat{c}}$  satisfies (A3) and we have $$\mathcal{D}(L)\subset \mathcal{F}\subset F.$$
In case (a), we see "$b$"-part in bilinear form as a perturbation to a symmetric Dirichlet form.

(iv)The semigroup can also be constructed by other methods. (see e.g. \cite{DR}, \cite{BDR}).

(v)  By (A3) we have that $\mathcal{E}$ is positivity preserving,  i.e.
$$\mathcal{E}(u,u^+)\geq 0\textrm{  }\forall u\in \mathcal{D}(L), $$
which can be obtained by the same arguments as in [36,  I  Proposition 4.4]. By (2.2) and (A3), we have for $u\in b\mathcal{D}(L), u\geq0$
$$\int Lud\mu=-\mathcal{E}(u,1)=-\int\langle Ab,\nabla u\rangle_H d\mu=-\int\langle Ab,\nabla (u+\varepsilon)\rangle_H d\mu \leq -\alpha \int (u+\varepsilon)d\mu .$$
Letting $\varepsilon\rightarrow0$, we have $\int Lud\mu\leq -\alpha \int ud\mu .$
$(P_t)_{t\in[0,T]}$ is a $C_0$-semigroup on $L^1(E;\mu)$.

(vi) All the conditions are satisfied by  the bilinear form considered in  [35, Section 4] and the operator in [13, Chapter II,III,IV] (see Section 6 below).

(vii) The notion of quasi-regularity for generalized Dirichlet forms analogously to \cite{MR} has been introduced in \cite{St}. By this and a technical assumption an associated $m$-tight special standard process can be constructed. We will use stochastic calculus associated with this process to conclude our probabilistic results (see Section 4 below).

(viii) Our assumptions (A1), (A2) are to make sure that the operator $L$ is associated to a bilinear form $\mathcal{E}$, which could be seen as a non-symmetric perturbation to the symmetric Dirichlet form $\mathcal{E}_A$. Under these assumptions we  prove the relation between generalized solutions (Definition 2.4), mild solutions and weak solutions in the sense of equation (2.7) in Proposition 2.7.  By this we obtain  a priori $L^\infty$-norm estimates for the solution of the nonlinear equation (3.1), which is essential to the proof of the existence of the solution to the nonlinear equation, since in our case the  condition (H4) on the nonlinearity $f$ is more general than all previous papers [10],[18],[26].

\vskip.10in

Let us recall the  notations  $\hat{F}, \mathcal{C}_T,\|\cdot\|_T$ associated with $\mathcal{E}^A$ from  \cite{BPS05}:
$\mathcal{C}_T:=C^1((0,T);L^2)\cap L^2(0,T;F)$, which turns out to be the appropriate space of test functions, i.e.
$$\aligned\mathcal{C}_T=\{&\varphi:[0,T]\times E\rightarrow \mathbb{R}|\varphi_t\in F \textrm{ for almost each } t, \int_0^T\mathcal{E}^A(\varphi_t,\varphi_t)dt<\infty,\\&t\rightarrow\varphi_t \textrm{ is differentiable in } L^2 \textrm {and } t\rightarrow\partial_t\varphi_t \textrm{ is } L^2-\textrm{continuous on }[0,T]\}.\endaligned$$
We also set $\mathcal{C}_{[a,b]}:=C^1([a,b];L^2)\cap L^2([a,b];F)$.
For $\varphi\in\mathcal{C}_T$, we define
$$\|\varphi\|_T:=(\sup_{t\leq T}\|\varphi_t\|_2^2+\int_0^T\mathcal{E}^A(\varphi_t)dt)^{1/2}.$$
$\hat{F}$ is the completion of $\mathcal{C}_T$ with respect to $\|\cdot\|_T$. By \cite{BPS05}, $\hat{F}=C([0,T];L^2)\cap L^2(0,T;F)$.
And for every $u\in\hat{F}$ there exists a sequence $u^n\in \mathcal{F}C_b^{\infty,T}, n\in \mathbb{N}$, such that $\int_0^T\mathcal{E}^A_1(u_t-u_t^n)dt\rightarrow0$. Here $$\mathcal{F}C_b^{\infty,T}:=\{f(t,l_1,...,l_m)| m\in \mathbb{N}, f\in C_b^\infty([0,T]\times \mathbb{R}^m), l_1,..., l_m\in E'\}.$$We also introduce the following space
$$W^{1,2}([0,T];L^2(E))=\{u\in L^2([0,T];L^2);\partial_t u\in L^2([0,T];L^2)\},$$
where $\partial_t u$ is the derivative of $u$ in the weak sense (see e.g. \cite{Ba10}).
\vskip.10in

\subsection{Linear Equations}

We consider the linear equation
\begin{equation}\aligned (\partial_t+L)u+f&=0,\qquad  0\leq t\leq T\\
u_T(x)&=\phi(x), \qquad x\in E\endaligned\end{equation}
where $f\in L^1([0,T];L^2(E,\mu)),\phi\in L^2(E,\mu)$.

By \cite{BPS05} we set $D_{A^{1/2}} \varphi:=A^{1/2}\nabla \varphi $ for any $\varphi\in \mathcal{F}C_b^\infty$, define $V_0=\{D_{A^{1/2}}\varphi:\varphi\in \mathcal{F}C_b^\infty\}$, and let $V$ be the closure of $V_0$ in $L^2(E;H,\mu)$.  And then we have the following results.
\vskip.10in

\th{Proposition 2.2} Assume (A1) holds.

(i) For every $u\in F$ there is a unique element of $V$, which we denote by $D_{A^{1/2}} u$ such that
$$\mathcal{E}^A(u,\varphi)=\int\langle D_{A^{1/2}} u(x),D_{A^{1/2}} \varphi(x)\rangle_H \mu(dx), \qquad \forall \varphi\in \mathcal{F}C_b^\infty.$$
One has $A^{1/2}A^{-1/2}D_{A^{1/2}} u(x)=D_{A^{1/2}} u(x)$. Moreover, the above formula extends to $u,v\in F$,
$$\mathcal{E}^A(u,v)=\int\langle D_{A^{1/2}} u(x),D_{A^{1/2}}v(x)\rangle_H \mu(dx). $$

(ii) Furthermore, if $u\in \hat{F}$, there exists a measurable function $\phi:[0,T]\times E\mapsto H$ such that $|{A^{1/2}}\phi |_H\in L^2([0,T]\times E)$ and $D_{A^{1/2}} u_t={A^{1/2}}\phi_t$ for almost all $t\in[0,T]$.

(iii)Let $u^n,u\in \hat{F}$ be such that $u^n\rightarrow u$ in $L^2((0,T)\times E)$ and $(D_{A^{1/2}} u^n)_n$ is a Cauchy-sequence in $L^2([0,T]\times E;H)$. Then $D_{A^{1/2}} u^n\rightarrow D_{A^{1/2}} u$ in $L^2((0,T)\times E;H)$, i.e. $D_{A^{1/2}}$ is closable as an operator from $\hat{F}$ into $L^2((0,T)\times E;H)$.

\proof See [4, Proposition 2.3].
$\hfill\Box$
\vskip.10in

For $u\in F,v\in bF$ we will denote
$$\mathcal{E}(u,v):=\int\langle D_{A^{1/2}} u(x),D_{A^{1/2}} v(x)\rangle_H \mu(dx)+\int \langle A^{1/2} b, D_{A^{1/2}} u\rangle_H v \mu(dx).$$

\th{Notation}By $\tilde{\nabla}u$ we denote  the set of all measurable functions $\phi:E\rightarrow H$, such that $A^{1/2}\phi=D_{A^{1/2}} u$ as elements of $L^2(\mu;H)$.

\subsection{Solution of the Linear Equation}

In this section we will introduce the concept of generalized solution and prove a generalized solution is equivalent to a mild solution in Proposition 2.7. Moreover, a generalized solution satisfies the weak relation (2.7). We don't use weak relation (2.7) as the definition of the solution. Since the solution is not in the domain of operator $L$, we can't choose it as a test function. It seems impossible to prove uniqueness of solution if we choose weak relation as the definition of the solution.

\th{Definition 2.3}[\emph{strong solution}] A function $u\in \hat{F}\cap L^1((0,T);\mathcal{D}(L))$ is called a strong solution of equation (2.3) with data $(\phi,f)$, if $t\mapsto u_t=u(t,\cdot)$ is $L^2$-differentiable on $[0,T],\partial _tu_t\in L^1((0,T);L^2)$ and the equalities in (2.3) hold in $L^2(\mu)$.

\th{Definition 2.4}[\emph{generalized solution}] A function $u\in \hat{F}$ is called a generalized solution of equation (2.3), if there exists a sequence of   $\{u^n\}$ consisting of  strong solutions with data $(\phi^n,f^n) $ such that $$\|u^n-u\|_T\rightarrow0,\|\phi^n-\phi\|_2\rightarrow0,\lim_{n\rightarrow\infty} f^n=f \textrm{ in } L^1([0,T];L^2(\mu)).$$

By (A3)and Remark 2.1 (v), for $0\leq t\leq T$,  $P_t$, as $C_0$-semigroup on $L^1(E;\mu)$,  can be restricted to a semigroup on $L^p(E;\mu)$ for all $p\in[1,\infty)$ by the Riesz-Thorin Interpolation Theorem and the restricted semigroup (denoted again by $P_t$ for simplicity) is strongly continuous on $L^p(E;\mu)$.

\th{Proposition 2.5} Assume that  (A1)-(A3) hold.

(i) Let $f\in C^1([0,T];L^p)$ for $p\in[1,\infty)$. Then
$w_t:=\int_t^TP_{s-t}f_sds\in C^1([0,T];L^p),$ and
$\partial_t w_t(x)=-P_{T-t}f_T(x)+\int_t^TP_{s-t}\partial_sf_s(x)ds.$

(ii) Assume that $\phi\in \mathcal{D}(L)$, $f\in C^1([0,T];L^2)$ and for each $t\in [0,T]$, $f_t\in \mathcal{D}(L)$ . Define
$u_t:=P_{T-t}\phi+\int_t^TP_{s-t}f_sds.$
Then $u$ is a strong solution of (2.3) and, moreover, $u\in C^1([0,T];L^2)$.

\proof By the same arguments as in [4, Proposition 2.6].  $\hfill\Box$

\vskip.10in
\th{Proposition 2.6} Assume that conditions (A1)-(A3) hold. If $f\in C^1([0,T], L^2(\mu))$ and $u$ is a strong solution for (2.3), it is a mild solution for (2.3) i.e.
$u_t=P_{T-t}\phi+\int_t^TP_{s-t}f_sds.$

\proof For fixed $t$, $\varphi\in \mathcal{D}(\hat{L})$
$(u_T,\hat{P}_{T-t}\varphi)-(u_t,\varphi)=\int_t^T(-Lu_s-f_s,\hat{P}_{s-t}\varphi)ds+
\int_t^T(u_s,\hat{L}\hat{P}_{s-t}\varphi)ds.$
Here $\hat{L}, \hat{P}_t$ denote the adjoints on $L^2(E,\mu)$ of $L, P_t$ respectively.
As $u$ is a strong solution, we deduce that
$(u_t,\varphi)=(P_{T-t}\phi+\int_t^TP_{s-t}f_sds,\varphi).$
Since $\mathcal{D}(\hat{L})$ is dense in $L^2$, the assertion follows.
$\hfill\Box$

\vskip.10in
\th{Proposition 2.7}Assume that conditions (A1)-(A3) hold, $f\in L^1([0,T];L^2)$ and $\phi\in L^2$. Then the equation (2.3) has a unique generalized solution  $u\in \hat{F}$
\begin{equation} u_t=P_{T-t}\phi+\int_t^TP_{s-t}f_sds.\end{equation}
The solution satisfies the three relations:
\begin{equation}\|u_t\|_2^2+2\int_t^T\mathcal{E}^A(u_s)ds\leq 2\int_t^T(f_s,u_s)ds+\|\phi\|_2^2+2\alpha\int_t^T\|u_s\|_2^2ds, \qquad 0\leq t\leq T,\end{equation}
\begin{equation}\|u\|_T^2\leq M_T(\|\phi\|_2^2+(\int_0^T\|f_t\|_2dt)^2),\end{equation}
\begin{equation}\int_0^T((u_t,\partial_t\varphi_t)+\mathcal{E}^A(u_t,\varphi_t)+\int\langle A^{1/2}b, D_{A^{1/2}}u_t
\rangle_H \varphi_td\mu)dt=\int_0^T(f_t,\varphi_t)dt+(\phi,\varphi_T)-(u_0,\varphi_0),\end{equation}
for any $\varphi\in b \mathcal{C}_T.$

Moreover, if $u\in\hat{F}$ is bounded and satisfies (2.7) for any $\varphi\in b \mathcal{C}_T$ with bounded $(f,\phi)$, then $u$ is a generalized solution given by (2.4). (2.7) can be extended easily to $\varphi\in bW^{1,2}([0,T];L^2)\cap L^2([0,T];F)$.

\proof Define $u$ by (2.4).
First assume that $\phi,f$ are bounded and satisfy the conditions of Proposition 2.5 (ii). Then, since $u$ is bounded and by Proposition 2.5 we know that $u$ is a strong solution of (2.3), hence it obviously satisfies (2.7). Furthermore, $u\in C^1([0,T];L^2)$.  Hence, actually $u\in b\mathcal{C}_T$ and consequently,
$$\int_{t_0}^T((u_t,\partial_tu_t)+\mathcal{E}^A(u_t,u_t)+\int\langle A^{1/2}b, D_{A^{1/2}}u_t
\rangle_H u_td\mu)dt=\int_{t_0}^T(f_t,u_t)dt+(\phi,u_T)-(u_{t_0},u_{t_0}).$$
By (2.2) we have $\int\langle A^{1/2}b, D_{A^{1/2}}u_t
\rangle_H u_td\mu\geq-\alpha\|u_t\|_2^2$  then we obtain
\begin{equation}\|u_t\|_2^2+2\int_t^T\mathcal{E}^A(u_s)ds\leq 2\int_t^T(f_s,u_s)ds+\|\phi\|_2^2+2\alpha\int_t^T\|u_s\|_2^2ds, \qquad 0\leq t\leq T.\end{equation}
As
$$\aligned\int_t^T(f_s,u_s)ds=&\int_t^T((f_s,P_{T-s}\phi)+(f_s,\int_s^TP_{r-s}f_rdr))ds
\\ \leq & M_0e^{T-t}(\|\phi\|_2\int_t^T\|f_s\|_2ds+\int_t^T(\|f_s\|_2\int_s^T\|f_r\|_2dr)ds),\endaligned$$
and
$\int_t^T\|u_s\|_2^2ds\leq M_{T-t}(\|\phi\|_2^2+(\int_0^T\|f_t\|_2dt)^2),$
 we obtain
$\|u_t\|_2^2+\int_t^T\mathcal{E}^{A}(u_s)ds\leq M_{T-t}(\|\phi\|_2^2+(\int_0^T\|f_t\|_2dt)^2).$
Hence,  it follows that
\begin{equation}\|u\|_T^2\leq M_T(\|\phi\|_2^2+(\int_0^T\|f_t\|_2dt)^2).\end{equation}
Here the constant $M_{T-t}$ may change from line to line, but it is independent of $f,\phi$. Now we will obtain the result for general data $\phi$ and $f$. Let $(f^n)_{n\in N}$ be a sequence of  functions in $bC^1([0,T];L^2(\mu))$ such that $f_t\in \mathcal{D}(L)$ for a.e. $t\in [0,T]$ and $\int_0^T\|f_t^n-f_t\|_2dt\rightarrow0$. (Such a sequence exists, since $\{\alpha_t g(x);\alpha_t\in C_0^\infty[0,T], g\in b\mathcal{D}(L)\}$ is dense in  $L^1([0,T];L^2)$). Take functions $(\phi^n)_{n\in N}\subset b\mathcal{D}(L)$ such that $\phi^n\rightarrow \phi$ in $L^2$. Let $u^n$ denote
the solution given by (2.4) with $f=f^n, \phi=\phi^n$.

By linearity, $u^n-u^m$ is associated with $(\phi^n-\phi^m, f^n-f^m)$. Since (2.9) implies that
$$\|u^n-u^m\|_T^2\leq M_T(\|\phi^n-\phi^m\|_2^2+(\int_0^T\|f_t^n-f_t^m\|_2dt)^2),$$
we  deduce that $(u^n)_{n\in N}$ is a Cauchy sequence in $\hat{F}$. Then $u=\lim_{n\rightarrow\infty}u^n$ in $\|\cdot\|_T$ is a generalized solution of (2.3) and (2.4) follows.

 Next we prove (2.5) (2.6) (2.7) for $u$.  For $\varphi\in b\mathcal{C}_T$, we have
 \begin{equation}\int_0^T((u_t^n,\partial_t\varphi_t)+\mathcal{E}^A(u_t^n,\varphi_t)+\int\langle A^{1/2}b, D_{A^{1/2}}u_t^n
\rangle_H  \varphi_td\mu)dt=\int_0^T(f_t^n,\varphi_t)dt+(\phi^n,\varphi_T)-(u_0^n,\varphi_0).\end{equation}
Since we have
$|\int_0^T\mathcal{E}^A(u_t^n-u_t,\varphi_t)dt|\leq (\int_0^T\mathcal{E}^A(u_t^n-u_t)dt)^{\frac{1}{2}}(\int_0^T\mathcal{E}^A(\varphi_t)dt)^{\frac{1}{2}}\rightarrow0,$
and
 $$\aligned |\int_0^T\int\langle A^{1/2}b, D_{A^{1/2}}(u_t^n-u_t)
\rangle_H  \varphi_td\mu dt|&\leq \|\varphi\|_\infty (\int_0^T\int|A^{1/2}b|^2_Hd\mu dt)^{\frac{1}{2}}(\int_0^T\int|D_{A^{1/2}} (u_t^n-u_t)|_H^2d\mu dt)^{\frac{1}{2}}
\\ &\rightarrow0,\endaligned$$
 we deduce (2.7) for any $\varphi\in b \mathcal{C}_T$.

Since $\|u_t^n\|_T\rightarrow\|u_t\|_T$, we conclude
$\lim_{n\rightarrow\infty}\int_0^T\mathcal{E}^A(u_t^n)dt=\int_0^T\mathcal{E}^A(u_t)dt.$
 As the relations (2.5), (2.6) hold for the approximating functions, by passing to the limit,  (2.5) and (2.6) follows for $u$.

[Uniqueness] Let $v\in \hat{F}$ be another generalized solution of (2.3) and let $(v^n)_{n\in N},(\tilde{\phi}^n)_{n\in N}, (\tilde{f}^n)_{n\in N}$ be the corresponding approximating sequences in the definition of the generalized solution.
By Proposition 2.8  $\sup_{t\in[0,T]}\|u^n_t-v^n_t\|_2^2\leq M_T(\|\phi^n-\tilde{\phi}^n\|_2^2+(\int_0^T\|f_t^n-\tilde{f}_t^n\|_2dt)^2).$
Letting $n\rightarrow\infty$, this implies $u=v$.

For the last result  we have $\forall t_0\geq 0,\varphi\in b\mathcal{C}_T$
\begin{equation}\int_{t_0}^T((u_t,\partial_t\varphi_t)+\mathcal{E}^{a,\hat{b}}(u_t,\varphi_t)+\int\langle b\sigma, D_\sigma u_t
\rangle \varphi_tdm)dt=\int_{t_0}^T(f_t,\varphi_t)dt+(\phi,\varphi_T)-(u_{t_0},\varphi_{t_0}). \end{equation}
For $t\geq \frac{1}{n} $,  define
$u_t^n:=n\int_0^{\frac{1}{n}}u_{t-s}ds,  f_t^n:=n\int_0^{\frac{1}{n}}f_{t-s}ds, \phi^n:=n\int_0^{\frac{1}{n}}u_{T-s}ds.$
It is easy to check that $u^n$ also fulfills (2.11) with $f^n,\phi^n$ i.e. for fixed $t_0\in (0,T],$ and for $n\geq \frac{1}{t_0}$,
$$\aligned &\int_{t_0}^T((u_t^n,\partial_t\varphi_t)+\mathcal{E}^{A}(u_t^n,\varphi_t)+\int\langle A^{1/2}b, D_{A^{1/2}} u_t^n\rangle \varphi_tdm)dt
=\int_{t_0}^T(f_t^n,\varphi_t)dt+(\phi^n,\varphi_T)-(u_{t_0}^n,\varphi_{t_0})dt.\endaligned$$
For the mild solution $v$ associated with $f,\phi$, the above relation also holds with $v^n$ replacing $u^n$. Hence we have
$$\int_{t_0}^T(((u-v)_t^n,\partial_t\varphi_t)+\mathcal{E}^{A}((u-v)_t^n,\varphi_t)+\int\langle A^{1/2}b, D_{A^{1/2}} (u-v)_t^n\rangle \varphi_tdm)dt=-((u-v)_{t_0}^n,\varphi_{t_0}).$$
Since $(u-v)_t^n\in b\mathcal{C}_{[\frac{1}{n},T]}$, the above equation holds with $(u-v)_t^n$ as a test function.
So we have $$\|(u-v)^n_{t_0}\|_2^2+2\int_{t_0}^T\mathcal{E}^{A}((u-v)_t^n,
(u-v)_t^n)dt\leq 2\alpha\int_{t_0}^T\|(u-v)^n_t\|_2^2dt.$$
By Gronwall's Lemma it follows that $\|(u-v)^n_{t_0}\|_2^2=0.$
Letting $n\rightarrow\infty$, we have $\|u_{t_0}-v_{t_0}\|_2=0$. Then letting $t_0\rightarrow0$, we have $\|u_0-v_0\|=0$. Therefore,
$u_t=P_{T-t}\phi+\int_t^TP_{s-t}f_sds$ is a generalized solution for (2.3).
$\hfill\Box$

\vskip.10in
We can prove the following basic relations for the linear equation which is essential to the following section. The basic idea of the proof comes from [4]. As the definition of the solution is different from [4] and our bilinear form is not symmetric, we need to apply some results in Proposition 2.7 and some  properties of the bilinear form $\mathcal{E}$ from (A1)-(A3) to conclude the following proposition. And here  we also use some new estimates to deal with the non-symmetric part  of the bilinear form $\mathcal{E}$. By Corollary A.4 and a modification of [4] we  prove Proposition 2.8. We omit it here. For more details, we refer to [40]. The proof of Corollary A.4 is included in Appendix A.
\vskip.10in

\th{Proposition 2.8} Let $u=(u^1,...,u^l)$ be a vector valued function where each component is a weak solution of the linear equation (2.3) associated to  data $f^i\in L^1([0,T];L^2),\phi^i\in L^2$    for $i=1,...,l$. By $\phi,f$ denote the vectors $\phi=(\phi^1,...,\phi^l),f=(f^1,...,f^l)$ and by $D_{A^{1/2}} u$ the matrix whose rows consist of $D_{A^{1/2}} u^i$. Then the following relations hold $\mu$-almost everywhere
\begin{equation}
 |u_t|^2+2\int_t^TP_{s-t}(|D_{A^{1/2}} u_s|_H^2)ds=P_{T-t}|\phi|^2+2\int_t^TP_{s-t}\langle u_s,f_s\rangle ds,
\end{equation}
and
\begin{equation}
 |u_t|\leq P_{T-t}|\phi|+\int_t^TP_{s-t}\langle \hat{u}_s,f_s\rangle ds.
\end{equation}
Here we write $\hat{x}=x/|x|$, for $x\in \mathbb{R}^l$, $x\neq0$ and $\hat{x}=0$, if $x=0$.

\section{The Non-linear Equation}
In the case of non-linear equations, we are going to treat systems of equations, with the unknown functions and their first-order derivatives mixed in the non-linear term of the equation. The non-linear term is a given measurable function $f:[0,T]\times E\times \mathbb{R}^l\times  H^l\rightarrow \mathbb{R}^l$, $l\in \mathbb{N}$.
We are going to treat the following system of equations.
\begin{equation}
 (\partial_t+L)u+f(\cdot,\cdot,u, D_{A^{1/2}} u)=0, \qquad u_T=\phi.
\end{equation}
The function $\phi$ is assumed to be in $L^2(E,d\mu;\mathbb{R}^l)$.
\vskip.10in
\th{Definition 3.1}[\emph{Generalized solution of the nonlinear equation}]
A generalized solution of equation (3.1) is a system $u=(u^1,u^2,...,u^l)$ of $l$ elements in $\hat{F}$, which has the property that each function $f^i(\cdot,\cdot, u, D_{A^{1/2}} u)$ belongs to $L^1([0,T];L^2(\mu))$ and such that
there is a sequence  $\{u_n\}$ which consists of strong solutions to (3.1) with data $(\phi_n,f_n) $ such that $$\|u_n-u\|_T\rightarrow0,\|\phi_n-\phi\|_2\rightarrow0\textrm{ and }\lim_{n\rightarrow\infty} f_n(\cdot,\cdot,u_n,D_{A^{1/2}} u_n)=f(\cdot,\cdot,u,D_{A^{1/2}} u) \textrm{ in } L^1([0,T];L^2(\mu)).$$
\vskip.10in
\th{Definition 3.2}[\emph{Mild solution}] A mild solution of equation (3.1) is a system $u=(u^1,u^2,...,u^l)$ of $l$ elements in $\hat{F}$, which has the property that each function $f^i(\cdot,\cdot, u, D_{A^{1/2}} u)$ belongs to $L^1([0,T];L^2(\mu))$ and such that
for every $i\in \{1,...,l\}$, the following equation holds
\begin{equation}u^i(t,x)=P_{T-t}\phi^i(x)+\int_t^TP_{s-t}f^i(s,\cdot,u_s,D_{A^{1/2}} u_s)(x)ds, \mu-a.e..\end{equation}
\vskip.10in
\th{Lemma 3.3} $u$ is a generalized solution of the nonlinear equation (3.1) if and only if it is a mild solution of equation (3.1).

\proof The assertion follows by Proposition 2.7.$\hfill\Box$
\vskip.10in
We will use the following notation
$|u|_H:=\sum|u^i|_H,$ for $ u\in L^2(E;H^l,d\mu),\|\phi\|_2^2:=\sum_{i=1}^l\|\phi^i\|_2^2,$ for $ \phi\in L^2(E,d\mu;\mathbb{R}^l),
\mathcal{E}(u,v):=\sum_{i=1}^l\mathcal{E}(u^i,v^i),\textrm{ }\mathcal{E}^A(u,v):=\sum_{i=1}^l\mathcal{E}^A(u^i,v^i),$ for $u,v\in F^l,$
$\|u\|_T^2:=\sup_{t\leq T}\|u_t\|_2^2+\int_0^T\mathcal{E}^A(u_t)dt,$ for $u\in \hat{F}^l.$
\subsection{The Case of Lipschitz Conditions}
In this subsection we consider a measurable function $f:[0,T]\times E\times \mathbb{R}^l\times H^l\rightarrow \mathbb{R}^l$ such that
\begin{equation}|f(t,x,y,z)-f(t,x,y',z')|\leq C(|y-y'|+|z-z'|_H),\end{equation}
with $t,x,y,y',z,z'$ arbitrary and $C$ a constant independent of $t,x$. We set $f^0(t,x):=f(t,x,0,0)$.
\vskip.10in
\th{Proposition 3.4} Assume that  conditions (A1)-(A3) hold and $f$ satisfies  condition (3.3), $f^0\in L^2([0,T]\times E,dt\times d\mu;\mathbb{R}^l)$ and $\phi\in L^2(E;\mathbb{R}^l)$. Then  equation (3.1) admits a unique solution $u\in\hat{F}^l$.
The solution satisfies the following estimate
$$\|u\|_T^2\leq e^{T(1+2C+C^2+2\alpha)}(\|\phi\|_2^2+\|f^0\|^2_{L^2([0,T]\times E)}).$$

\proof If $u\in \hat{F}^l$, then by relation (3.3) we have
$$\aligned |f(\cdot,\cdot,u,D_{A^{1/2}} u)|&\leq |f(\cdot,\cdot,u,D_{A^{1/2}} u)-f(\cdot,\cdot,0,0)|+|f(\cdot,\cdot,0,0)|
\\& \leq C(|u|+|D_{A^{1/2}} u|_H)+|f^0|.\endaligned$$
As $f^0 \in L^2([0,T]\times E,dt\times d\mu;\mathbb{R}^l)$ and $|D_{A^{1/2}} u|_H$ is an element of $L^2([0,T]\times E)$, we get $f(\cdot,\cdot,u,D_{A^{1/2}} u)\in  L^2([0,T]\times E;\mathbb{R}^l)$.

Now we define the operator $A:\hat{F}^l\rightarrow \hat{F}^l$ by
$(Au)^i_t=P_{T-t}\phi^i(x)+\int_t^T P_{s-t}f^i(s,\cdot,u_s,D_{A^{1/2}} u_s)(x)ds, i=1,...,l.$
Then  Proposition 2.7 implies that $Au\in \hat{F}^l$. In the following we write $f_{u,s}^i:=f^i(s,\cdot,u_s,D_{A^{1/2}} u_s).$ Since $(Au)^i_t-(Av)^i_t=\int_t^TP_{s-t}(f_{u,s}^i-f_{v,s}^i)ds$ is the mild solution with data $(f_u^i-f_v^i,0)$,  by the same argument as in Proposition 2.7 we have
$$\aligned&\|\int_t^TP_{s-t}(f_{u,s}^i-f_{v,s}^i)ds\|_{[t,T]}^2\leq M_T(\int_t^T\|f_{u,s}-f_{v,s}\|_2ds)^2
\\
\leq &M_T(T-t) \int_t^T(\|u_s-v_s\|_2^2+\||D_{A^{1/2}} u_s-D_{A^{1/2}} v_s|_H\|^2_2) ds
 \leq M_T(T-t) \|u-v\|_{[t,T]}^2,\endaligned$$
 where $M_T$ may change from line to line.
Here $\|u\|_{[T_a,T_b]}:=(\sup_{t\in[T_a,T_b]}\|u_t\|_2^2+\int_{T_a}^{T_b}\mathcal{E}^A(u_t)dt)^{\frac{1}{2}},$
where $0\leq T_a\leq T_b\leq T$. Fix $T_1$ sufficiently small such that $M_T(T-T_1)<1$. Then we have
:$$\|Au-Av\|_{[T_1,T]}^2<\|u-v\|^2_{[T_1,T]}.$$
Then there exists a unique $u_1\in \hat{F}_{[T_1,T]}$ such that $Au_1=u_1$ where $\hat{F}_{[T_a,T_b]}:=C([T_a,T_b];L^2)\cap L^2((T_a,T_b);F)$ for $T_a\in [0,T]$ and $T_b\in[T_a,T]$.
Then we can construct a solution on $[0,T]$ by iteration and uniqueness follows from the fixed point theorem.

In order to obtain the estimate in the statement, we write
$$\aligned |\int_t^T(f_{u,s},u_s)ds|
 \leq &\frac{1}{2}\int_t^T\|f_s^0\|_2^2ds+(\frac{1}{2}
+C+\frac{1}{2}C^2)\int_t^T\|u_s\|_2^2ds+\frac{1}{2}\int_t^T\mathcal{E}^A(u_s)ds.\endaligned$$
By relation (2.5) of Proposition 2.7 it follows that
$$\aligned\|u_t\|_2^2+2\int_t^T\mathcal{E}^A(u_s)ds&\leq2\int_t^T(f_{u,s},u_s)ds+\|\phi\|_2^2
+2\alpha\int_t^T\|u_s\|_2^2ds\\
&\leq\|\phi\|_2^2+\int_t^T\|f_s^0\|_2^2ds+(1+2C+C^2+2\alpha)\int_t^T\|u_s\|_2^2ds+\int_t^T\mathcal{E}^A(u_s)ds.
\endaligned$$
Now by Gronwall's lemma the desired estimate follows.

$\hfill\Box$
\vskip.10in
\subsection{The Case of Monotonicity Conditions}
Let $f:[0,T]\times E\times \mathbb{R}^l\times H^l\rightarrow \mathbb{R}^l$ be a measurable function and $\phi\in L^2(E,\mu;\mathbb{R}^l)$ be the final condition of (3.1). In this subsection  we impose the following conditions:
\vskip.10in
\no (H1) [\emph{Lipschitz condition in $z$}] There exists a fixed constant $C>0$ such that for $t,x,y,z,z'$ arbitrary
$|f(t,x,y,z)-f(t,x,y,z')|\leq C|z-z'|_H.$
\vskip.10in
\no(H2) [\emph{Monotonicity condition in $y$}]
For $x,y,y',z$ arbitrary, there exists a function $\mu\in L^1([0,T];\mathbb{R})$ such that
$\langle y-y',f(t,x,y,z)-f(t,x,y',z)\rangle
\leq \mu_t|y-y'|^2.$
We set $\alpha_t:=\int_0^t\mu_sds.$
\vskip.10in

\no (H3) [\emph{Continuity condition in $y$}]
For $t,x$ and $z$ fixed, the map
$\mathbb{R}^l\ni y\mapsto f(t,x,y,z)$ is continuous.

We need the following notation
$f^0(t,x):=f(t,x,0,0), f'(t,x,y):=f(t,x,y,0)-f(t,x,0,0),$
$f^{',r}(t,x):=\sup_{|y|\leq r}|f'(t,x,y)|.$
\vskip.10in
\no (H4) For each $r>0$, $f^{',r}\in L^1([0,T];L^2).$
\vskip.10in
\no (H5) $\|\phi\|_\infty<\infty,\|f^0\|_\infty<\infty.$

 As $\mu(E)<\infty$ we have $|\phi|\in L^2, |f^0|\in L^2([0,T];L^2)$. The conditions (H1), (H4), and (H5) imply that if $u\in \hat{F}$ is bounded, then $|f(u,D_{A^{1/2}} u)|\in L^1([0,T];L^2)$. Under the above conditions, even if $E$ is equal to a Hilbert space, it seems impossible to apply general monotonicity methods to the map $\mathcal{V}\ni u\mapsto f(t,\cdot,u(\cdot),D_{A^{1/2}}u)\in \mathcal{V}'$ because of lack of a suitable reflexive Banach space $\mathcal{V}$ such that $\mathcal{V}\subset\mathcal{H}\subset \mathcal{V}'$. Therefore, also here we proceed developing a hands-on approach to prove existence and uniqueness of solutions for equation (3.1) as done in \cite{BPS05}, \cite{S} and in particular, \cite{Z}.
\vskip.10in

\th{Lemma 3.5} In (H2) without loss of generality we  assume that $\mu_t\equiv0$.

\proof Let us make the change of variables $u_t^*:=\exp(\alpha_t)u_t$ and set $$\phi^*:=\exp(\alpha_T)\phi \qquad f_t^*(y,z):=\exp(\alpha_t) f_t(\exp(-\alpha_t)y, \exp(-\alpha_t)z)-\mu_t y$$ for the data. Next we can easily prove that $u$ is a solution associated to the data $(\phi,f)$ if and only if $u^*$ is a solution associated to the data $(\phi^*,f^*)$.
 It is obvious that (H1)-(H5) are satisfied. $\hfill\Box$
\vskip.10in

\th{Lemma 3.6} Assume that  conditions (A1)-(A3), (H1) and the following weaker form of condition (H2) (with $\mu_t\equiv0$) hold,
$$(H2')\langle y, f'(t,x,y)\rangle\leq 0,$$ for all $t,x,y$. If $u$ is a solution of (3.1), then there exists a constant $K$ which depends on $C, T,\alpha$ such that
$$\|u\|_T^2\leq K(\|\phi\|_2^2+\int_0^T\|f_t^0\|_2^2dt).$$
\proof Since $u$ is a solution of (3.1), we have by Proposition 2.7
$\|u_t\|_2^2+2\int_t^T\mathcal{E}^A(u_s)ds\leq 2\int_t^T(f_s,u_s)ds+\|u_T\|^2_2+2\alpha\int_t^T\|u_s\|_2^2ds.$
Conditions (H1) and (H2') yield
$$\aligned \langle f_s(u_s,D_{A^{1/2}} u_s),u_s\rangle=&\langle f_s(u_s,D_{A^{1/2}} u_s)-f_s(u_s,0)+f_s'(u_s)+f_s^0,u_s\rangle\\
\leq& (C|D_{A^{1/2}} u_s|_H+|f_s^0|)|u_s|.\endaligned$$
Hence, it follows
$$\aligned\|u_t\|_2^2+2\int_t^T\mathcal{E}^A(u_s)ds\leq &\int_t^T\mathcal{E}^A(u_s)ds+(C^2 +1+2\alpha)\int_t^T\|u_s\|_2^2ds+\int_t^T\|f_s^0\|_2^2ds+\|u_T\|^2_2.\endaligned.$$
Then by Gronwall's Lemma, the assertion follows.$\hfill\Box$
\vskip.10in
 By using Proposition 2.8, Lemma A.5 and Jesen's inequality, we have the following estimates. This can be done by a modification of the arguments in [4, Lemma 3.3] and we omit it here. For more details, we refer to [40].
\vskip.10in

\th{Lemma 3.7} Assume that the conditions (A1)-(A3), (H1) and (H2') hold. If $u$ is a solution of (3.1) , then  there exists a constant $K$, which depends on $C$ and $T$ such that
\begin{equation}
 \|u\|_\infty\leq K(\|\phi\|_\infty+\|f^0\|_\infty).
\end{equation}

\vskip.10in

Now we want to prove the following main theorem in this section.
The proof here is different from the finite dimension case since a unit ball in $H$ is not compact. And it is inspired from the probabilistic approach to prove the existence of the solution of the BSDE of \cite{BDHPS}.

\vskip.10in

\th{Theorem 3.8} Suppose the conditions (A1)-(A3), (H1)-(H5) hold. Then there exists a unique generalized solution of equation (3.1). And it satisfies the following estimates with constants $K_1$ and $K_2$ independent of $u,\phi,f$
$$\|u\|_T^2\leq K_1(\|\phi\|_2^2+\int_0^T\|f_t^0\|_2^2dt),$$
and
$$ \|u\|_\infty\leq K_2(\|\phi\|_\infty+\|f^0\|_\infty).$$
\proof [Uniqueness]

Let $u_1$ and $u_2$ be two solutions of equation (3.1). By using (2.5) for the difference $u_1-u_2$ we get
$$\aligned&\|u_{1,t}-u_{2,t}\|_2^2+2\int_t^T\mathcal{E}^A(u_{1,s}-u_{2,s})ds\\\leq& 2\int_t^T(f(s,\cdot,u_{1,s},D_{A^{1/2}} u_{1,s})-f(s,\cdot,u_{2,s},D_{A^{1/2}} u_{2,s}),u_{1,s}-u_{2,s})ds+2\alpha\int _t^T\|u_{1,s}-u_{2,s}\|_2^2ds\\
\leq &2\int_t^T C(|D_{A^{1/2}} u_{1,s}-D_{A^{1/2}} u_{2,s}|_H,|u_{1,s}-u_{2,s}|)ds+2\alpha\int _t^T\|u_{1,s}-u_{2,s}\|_2^2ds
\\ \leq & (C^2+2\alpha) \int _t^T\|u_{1,s}-u_{2,s}\|_2^2ds+\int_t^T\mathcal{E}^A(u_{1,s}-u_{2,s})ds.\endaligned$$
By Gronwall's lemma it follows that
$\|u_{1,t}-u_{2,t}\|_2^2=0,$
hence $u_1=u_2$.

[Existence] The existence will be proved in two steps.

\textbf{Step 1.} Suppose $f$ is bounded.
We define $M:=\sup|f(t,x,y,z)|.$

We need the following proposition.

\vskip.10in
\th{Proposition 3.9} If $f$ satisfies  the condition in \textbf{Step 1}, then for $v\in\hat{F}^l$, there exists a unique solution $u\in \hat{F}^l$ for the equation
$$(\partial_t+L)u+f(\cdot,\cdot,u, D_{A^{1/2}} v)=0, \qquad u_T=\phi.$$

Following the same arguments as in Lemma 3,5, we  assume that $2C^2+2\alpha+\mu_t\leq 0$,

For each $v\in\hat{F}^l$, we  define $Av=u$ where $u$ is the unique solution obtained by  Proposition 3.9. Let $v_1,v_2\in \hat{F}^l$. By
 applying (2.5) to the difference $u_1-u_2$ we get
$$\aligned&\|u_{1,t}-u_{2,t}\|_2^2+2\int_t^T\mathcal{E}^A(u_{1,s}-u_{2,s})ds\\\leq& 2\int_t^T(f(s,\cdot,u_{1,s},D_{A^{1/2}} v_{1,s})-f(s,\cdot,u_{2,s},D_{A^{1/2}} v_{2,s}),u_{1,s}-u_{2,s})ds+2\alpha\int _t^T\|u_{1,s}-u_{2,s}\|_2^2ds\\
\leq &2\int_t^T C(|D_{A^{1/2}} v_{1,s}-D_{A^{1/2}} v_{2,s}|_H,|u_{1,s}-u_{2,s}|)ds+\int _t^T(2\alpha+\mu_s)\|u_{1,s}-u_{2,s}\|_2^2ds
\\ \leq &\frac{1}{2}\int_t^T\mathcal{E}^A(v_{1,s}-v_{2,s})ds.\endaligned$$
Consequently we have $\|Av_1-Av_2\|_T\leq \frac{1}{2}\|v_1-v_2\|_T$. Then the fixed point $u$ of $A$ is the solution for (3.1).

\vskip.10in
\vspace{1mm}\noindent{\it Proof of Proposition 3.9}\quad We write $f(t,x,y)=f(t,x,y,D_{A^{1/2}} v)$.

We regularize $f$ with respect to the variable $y$ by convolution:
$$f_n(t,x,y,z)=n^l\int_{\mathbb{R}^l}f(t,x,y')\varphi(n(y-y'))dy'$$
where $\varphi$ is a smooth nonnegative function with support contained in the ball $\{|y|\leq 1\}$ such that $\int \varphi=1$.
Then $f=\lim_{n\rightarrow\infty}f_n$ and for each $n$, $\partial_{y_i}f_n$ are uniformly bounded.
Then each $f_n$ satisfies a Lipschitz condition with respect to both $y$ and $z$. Thus by Proposition 3.4 each $f_n$ determines  a solution $u_n\in \hat{F}^l$  of (3.1) with data $(\phi,f_n)$.
By the same arguments as in [34, Theorem 4.19], we have that each $f_n$ satisfies  conditions (H1) and (H2') with  $C=0$ and $\mu=0$. Since
$ |f_n(t,x,0,0)|
\leq n^l\int_{|y'|\leq \frac{1}{n}}|f(t,x,y')||\varphi(n(-y'))|dy'
\leq M,$
one deduces from Lemma 3.7 that $\|u_n\|_\infty\leq K$ and $\|u_n\|_T\leq K_T$.

Since the convolution operators approximate the identity uniformly on  compact sets, we get for fixed $t,x$,
$\lim_{n\rightarrow\infty}d'_{n,K}(t,x):=\sup_{|y|\leq K}|f(t,x,y)-f_n(t,x,y)|=0.$
Next we will show that $(u_n)_{n\in N}$ is a $\|\cdot\|_T$-Cauchy sequence. By (2.5)  for the difference $u_l-u_n$, we have
$$\aligned&\|u_{l,t}-u_{n,t}\|_2^2+2\int_t^T\mathcal{E}^A(u_{l,s}-u_{n,s})ds\\\leq& 2\int_t^T(f_l(s,\cdot,u_{l,s})-f_n(s,\cdot,u_{n,s}),u_{l,s}-u_{n,s})ds+2\alpha\int _t^T\|u_{l,s}-u_{n,s}\|_2^2ds
\\ \leq &2\int_t^T(|f_l(s,\cdot,u_{l,s})-f(s,\cdot,u_{l,s})|,|u_{l,s}-u_{n,s}|)ds+
2\int_t^T(|f_n(s,\cdot,u_{n,s})-f(s,\cdot,u_{n,s})|,|u_{l,s}-u_{n,s}|)ds\\&+
2\int_t^T(f(s,\cdot,u_{l,s})-f(s,\cdot,u_{n,s} ),u_{l,s}-u_{n,s})ds+2\alpha\int _t^T\|u_{l,s}-u_{n,s}\|_2^2ds\\\
\leq &2\int_t^T(d'_{l,K}(s,\cdot)+d '_{n,K}(s,\cdot),|u_{l,s}-u_{n,s}|)ds+2\alpha\int _t^T\|u_{l,s}-u_{n,s}\|_2^2ds\\
 \leq & \int_t^T\|d'_{l,K}(s,\cdot)\|_2^2ds+\int_t^T\|d '_{n,K}(s,\cdot)\|_2^2ds+(2+2\alpha)\int_t^T\|u_{l,s}-u_{n,s}\|_2^2ds,\endaligned$$
and that $\lim_{n\rightarrow\infty}\int_t^T\|d '_{n,r}(s,\cdot)\|_2^2ds=0$. Thus, for $l,n$ large enough, we get for an arbitrary $\varepsilon>0$
$$\|u_{l,t}-u_{n,t}\|_2^2+\int_t^T\mathcal{E}^A(u_{l,s}-u_{n,s})ds\leq \varepsilon+\tilde{K}\int_t^T\|u_{l,t}-u_{n,t}\|_2^2ds,$$
where $\tilde{K}$ depends on $C,M,\mu,\alpha$. It is easy to see that Gronwall's lemma then implies that $(u_n)_{n\in N}$ is a Cauchy-sequence in $\hat{F}$. Define $u:=\lim_{n\rightarrow\infty} u_n$ and take a subsequence $(n_k)_{k\in N}$ such that $u_{n_k}\rightarrow u$ a.e. We have $f(\cdot,\cdot,u_{n_k})\rightarrow f(\cdot,\cdot,u) \textrm{ in } L^2(dt\times d\mu).$
Since $\|u_{n_k}-u\|_T\rightarrow0$, we obtain
$\||D_{A^{1/2}} u-D_{A^{1/2}} u_{n_k}|_H\|_{L^2(dt\times d\mu)}\rightarrow0.$

We conclude
$$\aligned\lim_{k\rightarrow\infty}\|f_{n_k}(u_{n_k})-f(u)\|_{L^2(dt\times d\mu)}\leq &\lim_{k\rightarrow\infty}\|f_{n_k}(u_{n_k})-f(u_{n_k} )\|_{L^2(dt\times d\mu)}+\lim_{k\rightarrow\infty}\|f(u_{n_k})-f(u)\|_{L^2(dt\times d\mu)}\\\leq &\lim_{k\rightarrow\infty}\|d_{n_k,r}'\|_{L^2(dt\times d\mu)}+\lim_{k\rightarrow\infty}\|f(u_{n_k})-f(u)\|_{L^2(dt\times d\mu)}=0.\endaligned$$
By passing to the limit in the mild equation associated to $u_{n_k}$ with data $(\phi,f_{n_k})$, it follows that $u$ is the solution associated to $(\phi,f(u,D_{A^{1/2}} v))$.$\hfill\Box$

\vskip.10in
\textbf{Step 2.} Now we consider the general case. Let $r$ be a positive real number such that
$r\geq 1+K(\|\phi\|_\infty+\|f^0\|_\infty),$
where $K$ is the constant appearing in Lemma 3.7 (3.4).
Let $\theta_r$ be a smooth function such that $0\leq \theta_r\leq 1, \theta_r(y)=1$ for $|y|\leq r$ and $\theta_r(y)=0$ if $|y|\geq r+1$. For each $n\in N$, we set $q_n(z):=z\frac{n}{|z|_H\vee n}$ and
$h_n(t,x,y,z):=\theta_r(y)(f(t,x,y,q_n(z))-f_t^0)\frac{n}{f^{',r+1}\vee n}+f_t^0.$
We have $$\aligned |h_n(t,x,y,z)|&\leq |f(t,x,y,q_n(z))-f(t,x,y,0)+f(t,x,y,0)-f_t^0|1_{\{|y|\leq r+1\}}\frac{n}{f^{',r+1}\vee n}+f_t^0\\
&\leq C|q_n(z)|_H+\frac{nf^{',r+1}}{f^{',r+1}\vee n}+f_t^0\leq (1+C)n+f_t^0.\endaligned$$
We  easily show that $h_n$ satisfies (H1) and (H3). So, we only need to prove (H2).  For $y,y'\in \mathbb{R}^l$, if $|y|>r+1,|y'|>r+1$, the inequality is trivially satisfied and thus we concentrate on the case $|y'|\leq r+1$. We have
$$\aligned\langle y-y',h_n(t,x,y,z)-h_n(t,x,y',z)\rangle=&\theta_r(y)\frac{n}{f^{',r+1}\vee n}\langle y-y',f(t,x,y,q_n(z))-f(t,x,y',q_n(z))\rangle\\&+\frac{n}{f^{',r+1}\vee n}(\theta_r(y)-\theta_r(y'))\langle y-y',f(t,x,y',q_n(z))-f_t^0\rangle.\endaligned$$
The first term of the right hand side of the previous equality is negative. For the second term, we use that $\theta_r$ is $C(r)$-Lipschitz,
$$\aligned(\theta_r(y)-\theta_r(y'))\langle y-y',f(t,x,y',q_n(z))-f_t^0\rangle\leq & C(r)|y-y'|^2|f(t,x,y',q_n(z))-f_t^0|\\\leq& C(r)(Cn+f'_{r+1}(t))|y-y'|^2,\endaligned$$
and thus
$$\frac{n}{f^{',r+1}\vee n}(\theta_r(y)-\theta_r(y'))\langle y-y',f(t,x,y',q_n(z))-f_t^0\rangle\leq C(r)(C+1)n|y-y'|^2.$$
Then each $h_n$ satisfies  the assumptions in \textbf{Step 1}, and thus $u_n$ is the solution of (3.1) with data $(h_n,\phi)$. We have
$$\langle y, h_n'(t,x,y)\rangle=\langle y,h_n(t,x,y,0)-h_n(t,x,0,0)\rangle=\langle y, f(t,x,y,0)-f_t^0\rangle \frac{n\theta_r(y)}{f^{',r+1}\vee n}\leq0. $$
By Lemma 3.7, we also have $\|u_n\|_\infty\leq r-1$, $\|u_n\|_T\leq K_T$. So, $u_n$ is a solution with data $(f_n,\phi)$, where $f_n(t,x,y,z)=(f(t,x,y,q_n(z))-f_t^0)\frac{n}{f^{',r+1}\vee n}+f_t^0$. For this function (H2) is satisfied with $\mu_t=0$.
Conditions (H1) and (H2) yield
$$\aligned&|(f_l(u_l,D_{A^{1/2}} u_l)-f_n(u_n,D_{A^{1/2}} u_n),u_l-u_n)|\\
\leq &C(|D_{A^{1/2}} u_l-D_{A^{1/2}} u_n|_H,|u_l-u_n|)+|(f_l(u_n,D_{A^{1/2}} u_n)-f_n(u_n,D_{A^{1/2}} u_n),u_l-u_n)|.\endaligned$$

 For $n\leq l$, we have
 $$\aligned|f_l(u_n,D_{A^{1/2}} u_n)-f_n(u_n,D_{A^{1/2}} u_n)|\leq& 2C|D_{A^{1/2}} u_n|_H1_{\{|D_{A^{1/2}} u_n|_H\geq n\}}\\&+2C|D_{A^{1/2}} u_n|_H1_{\{f^{',r+1}> n\}}+2f^{',r+1}1_{\{f^{',r+1}> n\}}.\endaligned$$
Then we have
$$\aligned&\|u_{l,t}-u_{n,t}\|_2^2+2\int_t^T\mathcal{E}^A(u_{l,s}-u_{n,s})ds\\\leq& 2\int_t^T(f_l(u_{l,s},D_{A^{1/2}} u_{l,s})-f_n(u_{n,s},D_{A^{1/2}} u_{n,s}),u_{l,s}-u_{n,s})ds+2\alpha\int _t^T\|u_{l,s}-u_{n,s}\|_2^2ds\\
\leq &(C^2+2\alpha)\int_t^T\|u_l-u_n\|_2^2ds+\int_t^T\mathcal{E}^A(u_l-u_n)ds+8C(r-1)\int_t^T\int|D_{A^{1/2}} u_n|_H1_{\{|D_{A^{1/2}} u_n|_H\geq n\}}d\mu ds\\&+8C(r-1)\int_t^T|D_{A^{1/2}} u_n|_H1_{\{f^{',r+1}> n\}}d\mu ds+8C(r-1)\int_t^T\int f^{',r+1}1_{\{f^{',r+1}> n\}}d\mu ds.\endaligned$$
As $\|u_n\|_T^2\leq K_T$, we have $\int_0^T\||D_{A^{1/2}} u_n|_H\|_2^2ds\leq K_T$. Hence,
$$n^2\int_t^T\|1_{\{|D_{A^{1/2}} u_n|_H\geq n\}}\|_2^2ds\leq \int_t^T\||D_{A^{1/2}} u_n1_{\{|D_{A^{1/2}} u_n|_H\geq n\}}|_H\|_2^2ds\leq K_T.$$
As $\lim_{n\rightarrow\infty}\int_t^T\int_{\{f^{',r}>n\}}f^{',r}d\mu ds=0,$
and
$$\int_t^T\int_{\{f^{',r}>n\}}|D_{A^{1/2}} u_n|_Hd \mu dt\leq \|1_{\{f^{',r}>n\}}\|_{L^2(dt\times d\mu)}\||D_{A^{1/2}} u_n|_H\|_{L^2(dt\times d\mu)}\rightarrow0,$$
  for $n$ big enough we have
$$\|u_{l,t}-u_{n,t}\|_2^2+\int_t^T\mathcal{E}^A(u_{l,s}-u_{n,s})ds\\\leq (C^2+2\alpha)\int_t^T\|u_l-u_n\|_2^2ds+\varepsilon.$$
By Gronwalls' lemma it is easy to see that $(u_n)_{n\in N}$ is a Cauchy sequence in $\hat{F}^l$. Hence, $u:=\lim_{n\rightarrow\infty}u_n$ is well defined. We  find a subsequence such that $(u_{n_k},D_{A^{1/2}} u_{n_k})\rightarrow (u,D_{A^{1/2}} u)$ a.e.
$f(u_{n_k},D_{A^{1/2}} u)\rightarrow f(u,D_{A^{1/2}} u),$
and conclude that $$\aligned &|f_{n_k}(u_{n_k},D_{A^{1/2}} u_{n_k})-f(u,D_{A^{1/2}} u)|\\\leq&
1_{\{f^{',r}\leq n_k\}}|f(u,D_{A^{1/2}} u)-f(u_{n_k},q_{n_k}(D_{A^{1/2}} u_{n_k}))|\\&+1_{\{f^{',r}>n_k\}}[|f(u,D_{A^{1/2}} u)-f^0|+|f(u,D_{A^{1/2}} u)-f(u_{n_k},q_{n_k}(D_{A^{1/2}} u_{n_k}))|]\\ \leq &|f(u_{n_k},D_{A^{1/2}} u)-f(u_{n_k},q_{n_k}(D_{A^{1/2}} u_{n_k}))|+|f(u_{n_k},D_{A^{1/2}} u)-f(u,D_{A^{1/2}} u)|\\&+1_{\{f^{',r}>n_k\}}|f(u,D_{A^{1/2}} u)-f^0|\rightarrow 0 \textrm{ a.e.}.\endaligned$$
Since
$$\aligned &|f_{n_k}(u_{n_k},D_{A^{1/2}} u_{n_k})-f(u,D_{A^{1/2}} u)|\\\leq &|f(u,0)-f(u,D_{A^{1/2}} u)|+|f_{n_k}(u_{n_k},D_{A^{1/2}} u_{n_k})-f_{n_k}(u_{n_k},0)|+|f_{n_k}(u_{n_k},0)-f^0|+|f^0-f(u,0)|\\\leq &C(|D_{A^{1/2}} u|_H+|D_{A^{1/2}} u_{n_k}|_H)+2f^{',r},\endaligned$$
we have $$f_{n_k}(u_{n_k},D_{A^{1/2}} u_{n_k})\rightarrow f(u,D_{A^{1/2}} u)$$ in $L^1([0,T],L^2)$. We conclude $u$ is a solution of (3.1) associated to the data $(\phi,f)$.

$\hfill\Box$

\section{Martingale representation for the processes}
\subsection{Representation under $P^x$}
To relate the solution of non-linear equation (3.1) to backward stochastic differential equation, the most important part is to prove martingale representation theorem. In classical case, we could use the martingale representation theorem for Brownian motion. Now we want to extend this result for the process associated with the operator $L$. In \cite{QY}, they prove an abstract result about martingale representation theorem for Hunt process. Now we follow their idea to prove \emph{Fukushima representation property} mentioned below holds for our operator $L$ and extend their results to infinite dimensional case.

In order to obtain the results for the probabilistic part, we need that $\mathcal{E}$ is a generalized Dirichlet form in the sense of Remark 2.1 (iii) with $\hat{c}\equiv0$. There is a Markov process $X=(\Omega,\mathcal{F}_\infty,\mathcal{F}_t,X_t,P^x)$ which is properly associated in the resolvent sense with $\mathcal{E}$, i.e. $R_\alpha f:=E^x\int_0^\infty e^{-\alpha t}f(X_t)dt$ is $\mathcal{E}$-quasi-continuous $m$-version of the resolvent $G_\alpha$ of $\mathcal{E}$ for $\alpha>0$ and $f\in \mathcal{B}_b(E)\cap L^2(E;\mu)$. The coform $\hat{\mathcal{E}}$ introduced in Section 2 is a generalized Dirichlet form with the associated resolvent $(\hat{G}_\alpha)_{\alpha>0}$ and there exists an $m$-tight special standard process properly associated in the resolvent sense with $\hat{\mathcal{E}}$. We always assume that $(\mathcal{F}_t)_{t\geq0}$ is the (universally completed) natural filtration of $X_t$. From now on, we obtain all the results under the above assumptions.

As mentioned in Remark 2.1 (vii), such a process can be constructed by quasi-regularity ([36, IV. 1. Definition 1.7]) and a structural condition ([36, IV. 2. D3] on the domain $\mathcal{F}$ of the generalized Dirichlet form.

We will now introduce the spaces which will be relevant for our further investigations. Define
$$\aligned \mathcal{M}:=\{&M|M \textrm{ is a finite additive functional, } E^z[M_t^2]<\infty, E^z[M_t]=0\\& \textrm{ for }\mathcal{E}-q.e. z\in E \textrm{ and all } t\geq0\}.\endaligned$$
$M\in \mathcal{M}$ is called a \emph{martingale additive functional}(MAF). Let $M\in \mathcal{M}$. Then there exists an $\mathcal{E}$-exceptional set $N$, such that $(M_t,\mathcal{F}_t,P_z)_{t\geq 0}$ is a square integrable martingale for all $z\in E\backslash N$ (c.f. \cite{Tr 00}). Furthermore define
$$\dot{\mathcal{M}}=\{M\in\mathcal{M}|e(M)<\infty\}.$$
Here $e(M):=\frac{1}{2}\lim_{\alpha\rightarrow\infty}\alpha^2E^{\mu}[\int_0^\infty e^{-\alpha t}M_t^2dt]$ is the " energy" of $M$. The elements of $\dot{\mathcal{M}}$ are called martingale additive functional's (MAF) of finite energy. By [38, Theorem 2.10] $\dot{\mathcal{M}}$ is a real Hilbert space with inner product $e$.

Define for $k\in E$,
$$\mathcal{E}_k(u,v):=\int\frac{\partial u}{\partial k}\frac{\partial v}{\partial k}d\mu,\qquad u,v\in \mathcal{F}C_b^\infty.$$

$k\in E$ is called $\mu$-admissible if $(\mathcal{E}_k,\mathcal{F}C_b^\infty)$ is closable on $L^2(E;\mu)$.

We consider the following conditions:
\vskip.10in
(A4) There exists constants $c,C_1>0$ such that $c Id_H\leq A(z)\leq C_1 Id_H \textrm{ for all } z\in E.$ There exists  a countable dense subset  $\{e_k\}$ of $E'$, which is an orthonormal basis of $H$, consisting of $\mu$-admissible elements in $E$, and $u_k(\cdot):= _{E'}\langle e_k,\cdot\rangle_E\in \mathcal{F}$.
\vskip.10in

(A4') There exists  a countable dense subset  $\{e_k\}$ of $E'$, which is an orthonormal basis of $H$, consisting of $\mu$-admissible elements in $E$, and $u_k(\cdot)= _{E'}\langle e_k,\cdot\rangle_E\in \mathcal{F}$. Furthermore, $A(z)e_k=\lambda_k(z)e_k$ for some non-negative Borel measurable functions $\lambda_k$.
\vskip.10in

\th{Remark 4.1} Condition (A4) can be replaced by condition (A4') and all results below can be proved by the same argument. For simplicity, we will only give the proof under condition (A4).

By the existence of $\{e_k\}$, (A1) follows from [25, Proposition 3.8]. Set
$$\mathcal{F}C_b^\infty(\{e_k\}):=\{f(_{E'}\langle e_1,\cdot\rangle_E,...,_{E'}\langle e_m,\cdot\rangle_E)| m\in \mathbb{N}, f\in C_b^\infty(\mathbb{R}^m)\}.$$

\vskip.10in

(A5) The process $X$ associated with $\mathcal{E}$ above is a continuous conservative Hunt process in the state space $E\cup\{\partial\}$, $\alpha\hat{G}_\alpha$ is sub-Markovian and strongly continuous on $\mathcal{V}$, and $\hat{\mathcal{E}}$ is quasi-regular.
Furthermore, $\mathcal{F}C_b^\infty(\{e_k\})\subset\mathcal{F}$ and for $u\in \mathcal{F}$, there exists  a sequence $\{u_n\}\subset \mathcal{F}C_b^\infty(\{e_k\})$ such that $\mathcal{E}(u_n-u)\rightarrow0, n\rightarrow\infty$.

\vskip.10in
If $\mathcal{E}$ satisfies (A2) and (A4), we obtain$$\mathcal{E}(u,v):=\int \langle A(z)\nabla u(z),\nabla v(z)\rangle_H d\mu(z)+\int \langle A(z)b(z),\nabla u(z)\rangle_H v(z) d\mu(z), u\in F,v\in bF .$$
Again we set $D_{A^{1/2}}u:=A^{1/2}\nabla u$.

For  an initial distribution $\mu\in \mathcal{P}(E)$ ( where $\mathcal{P}(E)$ denotes all the probabilities on $E$, )  we will prove that the \emph{Fukushima reprensentation property} mentioned in \cite{QY} holds for $X$,  i.e.
there is an algebra $K(E)\subset\mathcal{B}_b(E)$ which generates the Borel $\sigma$-algebra $\mathcal{B}(E)$ and is invariant under $U^\alpha$ for $\alpha>0$, and there are countable continuous martingales $M^i,i\in \mathbb{N}$, over $(\Omega,\mathcal{F}^\mu,\mathcal{F}_t^\mu,P^\mu)$ such that for any potential $u=U^\alpha f$ where $\alpha>0$ and $f\in K(E)$, the martingale part $M^{[u]}$ of the semimartingale $u(X_t)-u(X_0)$ has a martingale representation in terms of $M^i$, that is, there are predictable processes $F_i,i\in \mathbb{N}$ on $(\Omega,\mathcal{F}^\mu,\mathcal{F}^\mu_t)$ such that $$M^{[u]}_t=\sum_{j=1}^\infty\int_0^tF_s^jdM_s^j\qquad P^\mu-a.e..$$

By [37, Theorem 4.5], if $\hat{G}_\alpha$ is sub-Markovian and strongly continuous on $\mathcal{V}$, Fukushima's decomposition holds for $u\in\mathcal{F}$. In this case we set $M^k:=M^{[u_k]}$, with $u_k(\cdot):=\langle e_k, \cdot\rangle_H.$ These martingales are called \emph{coordinate martingales}.

Let us first calculate the energy measure related to $\langle M^{[u]}\rangle, u\in \mathcal{F}C_b^\infty$. By  [38, formula (23)], for $g\in L^2(E,\mu)_b$, we have
$$\aligned\int \hat{G}_\gamma gd\mu_{\langle M^{[u]}\rangle}=&\lim_{\alpha\rightarrow\infty}\alpha(U_{\langle M^{[u]}\rangle}^{\alpha+\gamma}1,\hat{G}_\gamma g)\\=&\lim_{\alpha\rightarrow\infty}\lim_{t\rightarrow\infty}E_{\hat{G}_\gamma g\cdot \mu}(\alpha e^{-(\gamma+\alpha)t}\langle M^{[u]}\rangle_t)+\lim_{\alpha\rightarrow\infty}E_{\hat{G}_\gamma g\cdot \mu}(\int_0^\infty\langle M^{[u]}\rangle_t\alpha(\gamma+\alpha)e^{-(\gamma+\alpha)t}dt)\\=&\lim_{\alpha\rightarrow\infty}\lim_{t\rightarrow\infty}\alpha\langle \mu_{\langle M^{[u]}\rangle},e^{-(\gamma+\alpha)t}\int_0^t\hat{P}_s \hat{G}_\gamma gds\rangle\\&+ \lim_{\alpha\rightarrow\infty}\alpha(\gamma+\alpha)(\int_0^\infty e^{-(\gamma+\alpha)t}E_{\hat{G}_\gamma g\cdot \mu}((u(X_t)-u(X_0)-N_t^{[u]})^2)dt)\\=&\lim_{\alpha\rightarrow\infty}\alpha(\gamma+\alpha)(\int_0^\infty e^{-(\gamma+\alpha)t}E_{\hat{G}_\gamma g\cdot \mu}((u(X_t)-u(X_0))^2)dt)
\\=&\lim_{\alpha\rightarrow\infty}2\alpha(u-\alpha G_\alpha u, u\hat{G}_\gamma g)-\alpha(u^2, \hat{G}_\gamma g-\alpha\hat{G}_\alpha \hat{G}_\gamma g)
\\=&2(-Lu,u \hat{G}_\gamma g)-(-Lu^2,\hat{G}_\gamma g)
=2\mathcal{E}(u,u\hat{G}_\gamma g)-\mathcal{E}(u^2,\hat{G}_\gamma g)\\
=&2\mathcal{E}^A(u,u\hat{G}_\gamma g)-\mathcal{E}^A(u^2,\hat{G}_\gamma g)+2\int \langle Ab, \nabla u\rangle_H u\hat{G}_\gamma g \mu(dx)-\int \langle Ab, \nabla (u^2)\rangle_H \hat{G}_\gamma g \mu(dx)\\
=&2\mathcal{E}^A(u,u\hat{G}_\gamma g)-\mathcal{E}^A(u^2,\hat{G}_\gamma g)\\
=&2\int\langle A\nabla u,\nabla (u\hat{G}_\gamma g)\rangle_H d\mu-\int\langle A\nabla( u^2),\nabla (\hat{G}_\gamma g)\rangle_H d\mu\\
=&2\int\langle A \nabla u,\nabla u\rangle_H \hat{G}_\gamma gd\mu.\endaligned$$
Then by [38, Theorem 2.5], we have
$$\mu_{\langle M^{[u]}\rangle}=2\langle A \nabla u,\nabla u\rangle_H\cdot d\mu.$$

By [38, Proposition 2.19],  for $u\in\mathcal{F}C_b^\infty$ and $u=f(_{E'}\langle e_1,\cdot\rangle_E,...,_{E'}\langle e_m,\cdot\rangle_E)$, we have $$M^{[u]}_t=\sum_{i=1}^n\int_0^t\langle\nabla u(X_s),e_i\rangle_HdM_s^i.$$

Then by the same arguments as in [17, Theorem 3.1], we have that under $P^x$ for quasi every point $x$ ( where the exceptional set depends on $u,v$) , and every $u,v\in\mathcal{F}$,
\begin{equation} M^{[u]}_t=\sum_{i=1}^\infty\int_0^t\langle\nabla u(X_s),e_i\rangle_HdM_s^i.\end{equation}
Here $\sum_{i=1}^\infty\int_0^t\langle\nabla u(X_s),e_i\rangle_HdM_s^i=\lim_{n\rightarrow\infty}\sum_{i=1}^n\int_0^t\langle\nabla u(X_s),e_i\rangle_HdM_s^i$ in $(\dot{\mathcal{M}},e)$ and we have
\begin{equation}\langle M^{[u]},M^{[v]}\rangle_t=2\int_0^t\langle A(X_s)\nabla u(X_s),\nabla v(X_s)\rangle_H ds.\end{equation}
In particular,
 \begin{equation}\langle M^{i},M^{j}\rangle_t=2\int_0^ta_{ij}(X_s)ds,\end{equation}
 where $a_{ij}(z):=\langle A(z)e_i,e_j\rangle_H$.
\vskip.10in

 \th{Lemma 4.2}   Assume (A4)(A5) hold. For $u\in\mathcal{F}$ and $V_t$ is a continuous adapted process, with $|V_t|\leq M, \forall t,\omega$, we have for q.e. $x\in E$,
 \begin{equation}\int_0^tV_sdM_s^{[u]}=\sum_{i=1}^\infty\int_0^tV_s\langle\nabla u(X_s),e_i\rangle_HdM_s^i\qquad P^x-a.s..\end{equation}
 If (A4'), (A5) hold, then for some $\psi\in \tilde{\nabla}u$, we have
$$\int_0^tV_sdM_s^{[u]}=\sum_{i=1}^\infty\int_0^tV_s\langle\psi(X_s),e_i\rangle_HdM_s^i\qquad P^x-a.s..$$

 \proof By [38, Remark 2.2], for $\nu\in \hat{S}_{00}$ and  $B^n:=\sum_{i=1}^{n}\int_0^tV_s\langle\nabla u(X_s),e_i\rangle_HdM_s^i$ we have
 $$\aligned E^\nu(B^{n+m}-B^n)^2=& E^\nu(\sum_{i=n+1}^{n+m}\int_0^tV_s\langle\nabla u(X_s),e_i\rangle_HdM_s^i)^2
 \\=& E^\nu(\sum_{i,j=n+1}^{n+m}\int_0^tV_s^2a_{ij}(X_s)\langle\nabla u(X_s),e_i\rangle_H\langle\nabla u(X_s),e_j\rangle_Hds)\\\leq&C_1M^2E^\nu(\sum_{i=n+1}^{n+m}\int_0^t\langle\nabla u(X_s),e_i\rangle_H^2ds)\\\leq & C_1M^2e^t|U_1\nu|_\infty\sup_t\frac{1}{t}E^\mu(\sum_{i=n+1}^{n+m}\int_0^t\langle\nabla u(X_s),e_i\rangle_H^2ds)\\=&
  C_1M^2e^t|U_1\nu|_\infty\sum_{i=n+1}^{n+m}\int\langle\nabla u(z),e_i\rangle_H^2\mu(dz)\rightarrow0, \textrm{ as } n,m\rightarrow\infty.\endaligned$$
Then we define $\sum_{i=1}^\infty\int_0^tV_s\langle\nabla u(X_s),e_i\rangle_HdM_s^i:=\lim_{n\rightarrow\infty}B^n$ in $(\dot{\mathcal{M}},e)$. Furthermore, we have
  $$\aligned &E^\nu(\int_0^tV_sdM_s^{[u]}-\sum_{i=1}^n\int_0^tV_s\langle\nabla u(X_s),e_i\rangle_HdM_s^i)^2\\
  \leq&E^\nu\int_0^t\sum_{i,j=n+1}^\infty V_s^2a_{ij}(X_s)\langle\nabla u(X_s),e_i\rangle_H\langle\nabla u(X_s),e_j\rangle_Hds
  \\\leq &M^2C_1e^t|U_1\nu|_\infty \sup_t\frac{1}{t} E^\mu\int_0^t\sum_{i=n+1}^\infty \langle\nabla u(X_s),e_i\rangle_H^2ds\\\leq&
  M^2C_1e^t|U_1\nu|_\infty\int_0^t\sum_{i=n+1}^\infty \langle\nabla u(z),e_i\rangle_H^2\mu(dz)\rightarrow0, \textrm{ as } n\rightarrow\infty.\endaligned$$
So, we have $$\int_0^tV_sdM_s^{[u]}=\sum_{i=1}^\infty\int_0^tV_s\langle\nabla u(X_s),e_i\rangle_HdM_s^i\qquad P^\nu-a.s..$$
 Then by  [38, Theorem 2.5], the assertions follow.$\hfill\Box$
\vskip.10in

Moreover, by a modification of the proof of [31, Theorem 3.1] , we have the martingale representation theorem for $X$ which is similar to \cite{BPS05}.
\vskip.10in

\th{Theorem 4.3}Assume that (A4) or (A4') and (A5) hold. There exists some exceptional set $\mathcal{N}$ such that the following representation result holds: For  every bounded $\mathcal{F}_\infty$-measurable random variable $\xi$, there  exists predictable processes $\phi:[0,\infty)\times \Omega\rightarrow H$, such that for each probability measure $\nu$, supported by $E\setminus\mathcal{N}$, one has
$$\xi=E^\nu(\xi|\mathcal{F}_0)+\sum_{i=0}^\infty\int_0^\infty\langle\phi_s,e_i\rangle_HdM_s^{i}\qquad P^\nu-a.e.,$$
where $M^i=M^{[u_i]}$ with $u_i:=_{E'}\langle e_i, \cdot\rangle_{E}, i\in \mathbb{N}$ are the coordinate martingales, and
$$ E^\nu\int_0^\infty \langle A(X_s)\phi_s,\phi_s\rangle_Hds\leq \frac{1}{2}E^\nu\xi^2.$$
If another predictable process $\phi'$ satisfies the same relations under a certain measure $P^\nu$, then one has $A^{1/2}(X_t)\phi_t'=A^{1/2}(X_t)\phi_t, dt\times dP^\nu-a.e.$

\proof Suppose that $\mathcal{N}$ is some fixed exceptional set.  By $\mathcal{K}$ we denote the class of bounded random variables for which the statement holds outside this set.  We claim if $(\xi_n)\subset \mathcal{K}$ is a uniformly bounded increasing sequence and $\xi=\lim_{n\rightarrow\infty}\xi_n$, then $\xi\in\mathcal{K}$. Indeed, we have $E^x|\xi_n-\xi|^2\rightarrow0$. Let $\phi^{n}$ denote the process which represents $\xi_n$. Then
$$ E^x\int_0^\infty|\phi_s^{n}-\phi_s^{p}|_H^2ds\leq \frac{1}{c} E^x\int_0^\infty \langle A(X_s)(\phi_s^{n}-\phi_s^{p}),\phi_s^{n}-\phi_s^{p}\rangle_Hds\leq \frac{1}{2c}E^x|\xi_p-\xi_n|^2.$$
Now we want to pass to the limit with $\phi^n$ pointwise, so that the limit will become predictable. For each $l=0,1,...$ set
$$n_l(x):=\inf\{n|E^x(\xi-\xi_n)^2<\frac{1}{2^l}\},\quad\bar{\xi}_l:=\xi_{n_l(X_0)}.$$
Then one has $\bar{\xi}_l=\xi_{n_l(x)}$ on the set $\{X_0=x\}$, and $E^x(\xi-\bar{\xi}_l)^2<\frac{1}{2^l}$ for any $x\in\mathcal{N}^c$. The process which represents $\bar{\xi}_l$ is simply obtained by the formula $\bar{\phi}^l=\phi^{n_l(X_0)}$. Then define $\phi_s=\lim_{l\rightarrow\infty}\bar{\phi}^l_s$ in $H$. By the same argument as in Lemma 4.1, we have
$$\aligned E^x(\sum_{i=1}^\infty\int_0^\infty\langle\phi_s-\bar{\phi}_s^l,e_i\rangle_H dM_s^i)^2=&\lim_{k\rightarrow\infty}E^x(\sum_{i=1}^k\int_0^\infty\langle\phi_s-\bar{\phi}_s^l,e_i\rangle_H dM_s^i)^2\\
=&\lim_{k\rightarrow\infty}E^x(\sum_{i,j=1}^k\int_0^\infty a_{ij}(X_s)\langle\phi_s-\bar{\phi}_s^l,e_i\rangle_H\langle\phi_s-\bar{\phi}_s^l,e_j\rangle_H ds)\\
\leq&C_1E^x\int_0^\infty|\phi_s-\bar{\phi}_s^l|_H^2ds\rightarrow0.\endaligned$$
Therefore, we have $\xi\in\mathcal{K}.$

Let $K(E)\subset \mathcal{B}_b(E)$ be a countable set which is closed under  multiplication, generates the Borel
 $\sigma$-algebra $\mathcal{B}(E)$ and $U^\alpha (K(E))\subset K(E)$ for $\alpha \in \mathbb{Q}^+$. Such $K(E)$ can be constructed as follows.
 We choose a countable set $N_0\subset b\mathcal{B}(E)$ which generates the Borel $\sigma$-algebra $\mathcal{B}(E)$. Since $E$ as a separable Banach space is strongly Lindel\"{o}f, such a set $N_0$ can easily be constructed (see [25, Section 3.3]).  For $l\geq 1$ we define $N_{l+1}=\{g_1\cdot...\cdot g_k, U^\alpha g_1\cdot g_2\cdot...\cdot g_k,g_i\in N_l, k\in \mathbb{N}\cup \{0\},\alpha\in \mathbb{Q}^+\}$ and $K(E):=\cup_{l=0}^\infty N_l$ (c.f. [20 , Lemma 7.1.1]).

Let $\mathcal{C}_0$ be all $\xi=\xi_1\cdot\cdot\cdot\xi_n$ for some $n\in \mathbb{N}$,
$\xi_j=\int_0^\infty e^{-\alpha_jt}f_j(X_t)dt$, where $\alpha_j\in \mathbb{Q}^+$, $f_j\in K(E), j=1,...,n.$
Following the proof of [31, Lemma 2.2], we see that the universal completion of the $\sigma$-algebra generated by $\mathcal{C}_0$ is $\mathcal{F}_\infty$ .
By the first claim, a monotone class argument reduces the proof to the representation of a random variable in $\mathcal{C}_0$.

Let $\xi\in\mathcal{C}_0$. Following the same arguments as  in the proof of [31, Theorem 3.1] , we have
$N_t=E^x(\xi|\mathcal{F}_t)=\sum_m Z_t^m$ where the sum is  finite , and for each $m$, $Z^m=Z_t$ has the following form
$Z_t=V_tu(X_t)$
(the superscript $m$ will be dropped if no confusion may arise), where $V_t=\prod_{i=1}^{k'}\int_0^te^{-\beta_i s}g_i(X_s)ds$ and $u(x)=U^{\beta_1+...+\beta_k}(h_1(U^{\beta_2+...+\beta_k}h_2...(U^{\beta_k}h_k)...)$ for  $\beta_i\in \mathbb{Q}^+, g_i, h_i\in K(E)$. Obviously, $u\in K(E)$. Hence, by Fukushima's decomposition and Fukushima's representation property we have
\begin{equation}u(X_t)-u(X_0)=M_t^{[u]}+A_t^{[u]}=\sum_{j=1}^\infty\int_0^t \langle \nabla u(X_s), e_j\rangle_H dM_s^{j}+A_t^{[u]}\qquad P^x-a.s..\end{equation}
Then by the same arguments as in the proof of [31, Theorem 3.1]  and Lemma 4.2, we have
$$\aligned Z_t&=Z_0+\int_0^tu(X_s)dV_s+\int_0^tV_sdA_s^{[u]}+\int_0^tV_sdM_s^{[u]}\\
&=Z_0+\int_0^tu(X_s)dV_s+\int_0^tV_sdA_s^{[u]}+\sum_{i=1}^\infty\int_0^tV_s\langle\nabla u(X_s),e_i\rangle_HdM_s^i,\endaligned$$
$$N_t=\sum_{i=1}^\infty\int_0^tV_s\langle\nabla u(X_s),e_i\rangle_HdM_s^i \qquad P^x-a.s..$$
We  define $\phi_s=V_s\nabla u(X_s)$.
As (4.4) and (4.5) hold for every $x$ outside of an exceptional set of null capacity, the exceptional set $\mathcal{N}$ in the statement will be the union of all these exceptional sets corresponding to $u\in K(E)$ and the exceptional sets related to $V$ in Lemma 4.2.$\hfill\Box$
\vskip.10in

If in the preceding theorem, $\xi$ is nonnegative, we  drop the boundedness assumption.
\vskip.10in

\th{Corollary 4.4} Assume that (A4) or (A4') and (A5) hold. Let $\mathcal{N}$ be the set obtained in Theorem 4.3. For any $\mathcal{F}_\infty$-measurable nonnegative random variable $\xi\geq0$ there exists a predictable process $\phi:[0,\infty)\times \Omega\rightarrow H$ such that
$$\xi=E^x(\xi|\mathcal{F}_0)+\sum_{i=0}^\infty\int_0^\infty\langle\phi_s,e_i\rangle_HdM_s^{i}\qquad P^x-a.e.,$$
where $M^i, i\in \mathbb{N}$, are as in Theorem 4.3, and
$$E^\nu\int_0^\infty\langle A(X_s)\phi_s,\phi_s\rangle_Hds\leq \frac{1}{2}E^\nu\xi^2,$$
for each point $x\in\mathcal{N}^c$ such that $E^x\xi<\infty$.

If another predictable process $\phi'$ satisfies the same relations under a certain measure $P^x$, then one has $A^{1/2}(X_t)\phi_t'=A^{1/2}(X_t)\phi_t, dt\times dP^x-a.e.$

\vskip.10in

\subsection{Representation under $P^\mu$}
As usual we set $\int_0^t\psi_s.dM_s=\sum_{i=0}^\infty\int_0^t\langle\psi_s,e_i\rangle_HdM_s^{i}$.
\vskip.10in

\th{Lemma 4.5} Assume that (A1)-(A3), (A5) and (A4) or (A4') hold. If $u\in\mathcal{D}(L)$, $\psi\in \tilde{\nabla} u$ , then
$$u(X_t)-u(X_0)=\int_0^t\psi(X_s).dM_s+\int_0^tLu(X_s)ds\qquad P^\mu- a.s..$$
\proof Corollary 4.4 and (4.1) imply the assertion.$\hfill\Box$
\vskip.10in

The aim of the rest of this section is to extend this representation to time dependent functions $u(t,x)$.
\vskip.10in

\th{Lemma 4.6} Assume that (A1)-(A3), (A5) and (A4) or (A4') hold. Let $u:[0,T]\times E\rightarrow \mathbb{R}$ be such that

(i) $\forall s,u_s\in \mathcal{D}(L)$ and $s\rightarrow Lu_s$ is continuous in $L^2$.

(ii) $u\in C^1([0,T];L^2)$.

Then clearly $u\in\mathcal{C}_T$, and, moreover,  for any $\psi\in \tilde{\nabla }u$ and any $s,t>0$ such that $s+t<T$, the following relation holds $P^\mu$-a.s.
$$u(s+t,X_t)-u(s,X_0)=\int_0^t\psi_{s+r}(X_r).dM_r+\int_0^t(\partial_s+L)u_{s+r}(X_r)dr.$$
\proof We prove the above relation with $s=0$, the general case being similar. Let $0=t_0<t_1<...<t_p=t$ be a partition of the interval $[0,t]$ and write
$u(t,X_t)-u(0,X_0)=\sum_{n=0}^{p-1}(u(t_{n+1},X_{t_{n+1}})-u(t_n,X_{t_n})).$
Then, on account of the preceding lemma,
$$\aligned u(t_{n+1},X_{t_{n+1}})-u(t_n,X_{t_n})=&u(t_{n+1},X_{t_{n+1}})-u(t_{n+1},X_{t_{n}})+u(t_{n+1},X_{t_n})-u(t_n,X_{t_n})\\
=&\int_{t_n}^{t_{n+1}}\nabla u_{t_{n+1}}(X_s).dM_s+\int_{t_n}^{t_{n+1}}Lu_{t_{n+1}}(X_s)ds+\int_{t_n}^{t_{n+1}}\partial_su_s(X_{t_n})ds,\endaligned$$
where the last integral is obtained by using the Leibnitz-Newton formula for the $L^2$-valued function $s\rightarrow u_s$. Further we estimate the $L^2$-norm of the differences between each term in the last expression and the similar terms corresponding to the formula we have to prove. Now we use $\mu P_t\leq \mu$,  i.e. $\int P_tfd\mu\leq \int fd\mu$ for $f\in\mathcal{B}^+$ since $\alpha\hat{G}_\alpha$ is sub-Markovian, to obtain
$$\aligned &E^\mu(\int_{t_n}^{t_{n+1}}\nabla u_{t_{n+1}}(X_s).dM_s-\int_{t_n}^{t_{n+1}}\nabla u_{s}(X_s).dM_s)^2\\
=&E^\mu\int_{t_n}^{t_{n+1}}\langle A(X_s)(\nabla u_{t_{n+1}}(X_s)-\nabla u_{s}(X_s)),\nabla u_{t_{n+1}}(X_s)-\nabla u_{s}(X_s)\rangle_Hds
\leq\int_{t_n}^{t_{n+1}}\mathcal{E}^A(u_{t_{n+1}}-u_s)ds\endaligned$$
Since $s\rightarrow Lu_s$ is continuous in $L^2$, it follows that $s\rightarrow u_s$ is continuous w.r.t. $\mathcal{E}_1^A$-norm. Hence the integrand of the last integral  is uniformly small, provided the partition is fine enough. From this we deduce that
$\sum_{n=0}^{p-1}\int_{t_n}^{t_{n+1}}\nabla u_{t_{n+1}}(X_s).dM_s\rightarrow \int_0^t\nabla u_{r+s}(X_r).dM_r,$
as the mesh of the partitions tends to zero.
By Minkowski's inequality, we obtain
$$\aligned &(E^\mu(\sum_{n=0}^{p-1}\int_{t_n}^{t_{n+1}}(Lu_{t_{n+1}}-Lu_s)(X_s)ds)^2)^{1/2}\\
\leq&\sum_{n=0}^{p-1}\int_{t_n}^{t_{n+1}}(E^\mu(Lu_{t_{n+1}}-Lu_s)^2(X_s))^{1/2}ds
\leq\sum_{n=0}^{p-1}\int_{t_n}^{t_{n+1}}\|Lu_{t_{n+1}}-Lu_s\|_2ds,\endaligned$$
where the integrand again is a uniformly small quantity.

And finally
$$\aligned&(E^\mu(\sum_{n=0}^{p-1}\int_{t_n}^{t_{n+1}}(\partial_su_s(X_{t_n})-\partial_su_s(X_s))ds)^2)^{1/2}
\leq \sum_{n=0}^{p-1}\int_{t_n}^{t_{n+1}}(E^\mu(\partial_su_s(X_{t_n})-\partial_su_s(X_s))^2)^{1/2}ds\\
=&\sum_{n=0}^{p-1}\int_{t_n}^{t_{n+1}}(E^\mu(\partial_su_s(X_{t_n})^2+P_{s-t_n}(\partial_su_s)^2(X_{t_n})-2\partial_su_s(X_{t_n})(P_{s-t_n}
\partial_su)(X_{t_n})))^{1/2}ds\\
=&\sum_{n=0}^{p-1}\int_{t_n}^{t_{n+1}}(E^\mu((\partial_su_s(X_{t_n})-(P_{s-t_n}\partial_su_s)(X_{t_n}))^2+(P_{s-t_n}(\partial_su_s)^2(X_{t_n})-
((P_{s-t_n}\partial_su_s)(X_{t_n}))^2)))^{1/2}ds\\
\leq&\sum_{n=0}^{p-1}(\int_{t_n}^{t_{n+1}}\int(\partial_su_s
-P_{s-t_n}\partial_su_s)^2+P_{s-t_n}(\partial_su_s)^2-(P_{s-t_n}\partial_su_s)^2d\mu)^{1/2}ds.\endaligned$$
From the hypotheses it follows that this  tends also to zero as the mesh of the partitions goes to zero. Hence the assertion follows.$\hfill\Box$
\vskip.10in

\th{Theorem 4.7} Assume that (A1)-(A3), (A5) and (A4) or (A4') hold. Let $f\in L^1([0,T];L^2)$ and $\phi\in L^2(E)$ and define
$$u_t:=P_{T-t}\phi+\int_t^TP_{s-t}f_sds.$$
Then for each $\psi\in\tilde{\nabla}u$ and for each  $s\in[0,T]$, the following relation holds $P^\mu$-a.s.
$$u(s+t,X_t)-u(s,X_0)=\int_0^t\psi(s+r,X_r).dM_r-\int_0^tf(s+r,X_r)dr.$$ Furthermore,  if $u$ is a generalized solution of PDE (3.1), the following BSDE holds $P^\mu$-a.s.
$$\aligned u(t,X_{t-s})=&\phi(X_{T-s})+\int_t^Tf(r,X_{r-s},u(r,X_{r-s}),A^{1/2}\nabla u (r,X_{r-s}))dr\\&-\int_{t-s}^{T-s}\psi(s+r,X_r).dM_r.\endaligned$$
\proof  First assume that $\phi$ and $f$ satisfy the conditions in Proposition 2.5 (ii). Then we have that $u$ satisfies the condition in Lemma 4.6 and  by  Lemma 4.6, the assertion follows. For the general case we choose $u^n$ associated with $(f^n,\phi^n)$ as in Proposition 2.7. Then  $\|u^n-u\|_T\rightarrow0$ as $n\rightarrow\infty$.
For $u^n$ we have
\begin{equation}u^n(s+t,X_t)-u^n(s,X_0)=\int_0^t\nabla u^n_{s+r}(X_r).dM_r-\int_0^tf^n(s+r,X_r)dr.\end{equation}
As
$$\aligned &E^\mu|\int_0^t(\nabla u^n_{s+r}(X_r)-\nabla u^p_{s+r}(X_r)).dM_r|^2\\\leq& E^\mu\int_0^t\langle A(X_r)(\nabla u^n_{s+r}(X_r)-\nabla u^p_{s+r}(X_r)),\nabla u^n_{s+r}(X_r)-\nabla u^p_{s+r}(X_r)\rangle_Hdr
\leq \int_0^t\mathcal{E}^A(u^n_{s+r}-u^p_{s+r})dr,\endaligned$$
letting $n\rightarrow\infty$ in (4.6), we obtain the assertions.
$\hfill\Box$

\vskip.10in
\section{BSDE's and Weak Solutions}
The set $\mathcal{N}$ obtained in Theorem 4.3 will be fixed throughout this section. By Theorem 4.3, we  solve the BSDE under all measures $P^x$, $x\in\mathcal{N}^c$,
 at the same time. We will treat systems of $l$ equations, $l\in \mathbb{N}$, associated to $\mathbb{R}^l$-valued functions $f:[0,T]\times \Omega\times \mathbb{R}^l\times H^l\mapsto \mathbb{R}^l$, assumed to be predictable.
 This means that we consider the map $(s,\omega)\mapsto f(s,\omega,\cdot,\cdot)$ as a process which is predictable with respect to the canonical filtration of our process $(\mathcal{F}_t)$.

\vskip.10in
\th{Lemma 5.1} Assume that (A4) or (A4') and (A5) hold. Let $\xi$ be an $\mathcal{F}_T$-measurable random variable and $f:[0,T]\times \Omega\mapsto \mathbb{R}$ an $(\mathcal{F}_t)_{t\geq0}$-predictable process. Let $D$ be the set of all points $x\in \mathcal{N}^c$ for which the following integrability condition holds
$$E^x(|\xi|+\int_0^T|f(s,\omega)|ds)^2<\infty.$$
Then there exists a pair $(Y_t,Z_t)_{0\leq t\leq T}$ of predictable processes $Y:[0,T)\times \Omega\mapsto \mathbb{R}, Z:[0,T)\times \Omega\mapsto H$, such that under all measures $P^x$, $x\in D$, they have the following properties:

\no (i) $Y$ is continuous;

\no  (ii) $Z$ satisfies the integrability condition
$\int_0^T|A^{1/2}(X_t)Z_t|_H^2dt<\infty,  P^x-a.s.;$
(iii) the local martingale obtained integrating $Z$ against the coordinate martingales, i.e. $\int_0^tZ_s.dM_s$, is a uniformly integrable martingale;

\no (iv)
$Y_t=\xi+\int_t^Tf(s,\omega)ds-\int_t^TZ_s.dM_s, P^x-a.s., 0\leq t\leq T.$
If another pair $(Y_t',Z_t')$ of predictable processes satisfies the above conditions (i),(ii),(iii),(iv), under a certain measure $P^\nu$ with the initial distribution $\nu$ supported by $D$, then one has $Y.=Y.', P^\nu-a.s.$ and $A^{1/2}(X_t)Z_t=A^{1/2}(X_t)Z_t', dt\times P^\nu-a.s..$

\proof The representations of the positive and negative parts of the random variable $\xi+\int_0^Tf_sds$ give us a predictable process $Z$ such that
$\xi+\int_0^Tf_sds=E^{X_0}(\xi+\int_0^Tf_sds)+\int_0^TZ_s.dM_s.$
Then we obtain the desired process $Y$ by the formula
$Y_t=E^{X_0}(\xi+\int_0^Tf_sds)+\int_0^TZ_s.dM_s-\int_0^tf_sds.$
$\hfill\Box$

\vskip.10in
\th{Definition 5.2} Let $\xi$ be an $\mathbb{R}^l$-valued, $\mathcal{F}_T$-measurable, random variable and $f:[0,T]\times \Omega\times \mathbb{R}^l\times H^l\mapsto \mathbb{R}^l$ a measurable $\mathbb{R}^l$-valued
 function such that $(s,\omega)\mapsto f(s,\omega,\cdot,\cdot)$  is a predictable  process. Let $p>1$ and $\nu$ be a probability measure supported by $\mathcal{N}^c$ such that $E^\nu|\xi|^p<\infty$.
We say that a pair $(Y_t,Z_t)_{0\leq t\leq T}$ of predictable processes $Y:[0,T)\times \Omega\mapsto \mathbb{R}^l$, $Z:[0,T)\times\Omega\mapsto  H^l$
 is a solution of the BSDE in $L^p(P^\nu)$ with data $(\xi,f)$ provided $Y$ is continuous under $P^\nu$ and it satisfies both the integrability conditions
$\int_0^T|f(t,\cdot,Y_t,{A^{1/2}}(X_t)Z_t)|dt<\infty, P^\nu-a.s.,$
and
$E^\nu(\int_0^T|{A^{1/2}}(X_t)Z_t|_H^2dt)^{p/2}<\infty,$
and the following equation holds
\begin{equation}Y_t=\xi+\int_t^Tf(s,\omega, Y_s,{A^{1/2}}(X_s)Z_s)ds-\int_t^TZ_s.dM_s,\qquad P^\nu-a.s., 0\leq t\leq T.\end{equation}

\vskip.10in

Now to prove the existence of solution of BSDE we give the condition on nonlinear term $f$.
Let $f:[0,T]\times \Omega\times \mathbb{R}^l\times H^l\mapsto \mathbb{R}^l$ be a measurable $\mathbb{R}^l$-valued function such that
 $(s,\omega)\mapsto f(s,\omega,\cdot,\cdot)$ is predictable and it satisfies the following conditions:

\no($\Omega1$) [Lipschitz condition in $z$]
There exists a  constant $C>0$ such that for all $t,\omega,y,z,z'$
$|f(t,\omega,y,z)-f(t,\omega,y,z')|\leq C|z-z'|_H.$

\no($\Omega2$) [Monotonicity condition in $y$]
For $\omega,y,y',z$ arbitrary, there exists a function $\mu_t\in L^1([0,T],\mathbb{R})$ such that
$\langle y-y',f(t,\omega,y,z)-f(t,\omega,y',z)\rangle
\leq \mu_t|y-y'|^2,$
and  $\alpha_t:=\int_0^t\mu_sds.$

\no($\Omega3$) [Continuity condition in $y$]
For $t,\omega$ and $z$ fixed, the map
$y\mapsto f(t,\omega,y,z)$ is continuous.

We need the following notation
$f^0(t,\omega):=f(t,\omega,0,0), f'(t,\omega,y):=f(t,\omega,y,0)-f(t,\omega,0,0),$ $f^{',r}(t,\omega):=\sup_{|y|\leq r}|f'(t,\omega,y)|.$

Let $\xi$ be an $\mathbb{R}^l$-valued, $\mathcal{F}_T$-measurable, random variable and, for each $p>0$ let $A_p$ denote the set of all points $x\in\mathcal{N}^c$ such that
\begin{equation}E^x\int_0^Tf_t^{',r}dt<\infty, \qquad \forall r\geq 0,\end{equation}
and
$E^x(|\xi|^p+(\int_0^T|f^0(s,\omega)|ds)^p)<\infty.$
Let $A_\infty$ denote the set of points $x\in\mathcal{N}^c$ for which (5.2) holds and with the property that $|\xi|,|f^0|\in L^\infty(P^x)$.

The method to prove the following proposition is standard and is a modification of \cite{BDHPS}, we omit it here.

\vskip.10in

\th{Proposition 5.3} Assume that (A4) or (A4') and (A5) holds. Under conditions $(\Omega1),(\Omega2),(\Omega3)$ there exists a pair $(Y_t,Z_t)_{0\leq t\leq T}$ of predictable processes $Y:[0,T)\times \Omega\mapsto \mathbb{R}^l, Z:[0,T)\times\Omega\mapsto \mathbb{R}^l\otimes \mathbb{R}^d$ that forms a solution of the BSDE (5.1) in $L^p(P^x)$ with data $(\xi,f)$ for each point $x\in A_p$. Moreover, the following estimate holds with some constant $K$ that depends only on $c,C,\mu$ and $T$,
$$E^x(\sup_{t\in[0,T]}|Y_t|^p+(\int_0^T|{A^{1/2}}(X_t)Z_t|_H^2dt)^{p/2})\leq KE^x(|\xi|^p+(\int_0^T|f^0(s,\omega)|ds)^p), \qquad x\in A_p.$$
If $x\in A_\infty$, then $\sup_{t\in[0,T]}|Y_t|\in L^\infty(P^x)$.

If $(Y_t',Z_t')$ is another solution in $L^p(P^x)$ for some point $x\in A_p$, then one has $Y_t=Y_t'$ and ${A^{1/2}}(X_t)Z_t={A^{1/2}}(X_t)Z_t',dt\times P^x-a.s.$.

\vskip.10in

We shall now look at the connection between the solutions of BSDE's introduced in this section and the PDE's studied in Section 3. In order to do this we have to consider BSDE's over time intervals like $[s,T]$, with $0\leq s\leq T$. Since the present approach is based on the theory of  time homogeneous Markov processes,  we have to discuss solutions over the interval $[s,T]$, where the process and the coordinate martingales are indexed by a parameter in the interval $[0,T-s]$.

Let us give a formal definition for the natural notion of solution over a time interval $[s,T]$. Let $\xi$ be an $\mathcal{F}_{T-s}$-measurable, $\mathbb{R}^l$-valued, random variable and $f:[s,T]\times \Omega\times \mathbb{R}^l\times H^l\rightarrow \mathbb{R}^l$ an $\mathbb{R}^l$-valued, measurable map such that $(f(s+l,\omega,\cdot,\cdot))_{l\in[0,T-s]}$ is predictable with respect to $(\mathcal{F}_l)_{l\in[0,T-s]}$. Let $\nu$ be a probability measure supported by $\mathcal{N}^c$ such that $E^\nu|\xi|^p<\infty$. We say a pair $(Y_t,Z_t)_{s\leq t\leq T}$ of processes $Y:[s,T]\times\Omega\rightarrow \mathbb{R}^l, Z:[s,T]\times\Omega\rightarrow H^l$ is a solution in $L^p(P^\nu)$ of the BSDE over the interval $[s,T]$ with data $(\xi,f)$, provided  they have the property that reindexed as $(Y_{s+l},Z_{s+l})_{l\in[0,T-s]}$ these processes are $(\mathcal{F}_l)_{l\in[0,T-s]}$-predictable, $Y$ is continuous and together they satisfy the integrability conditions
$\int_s^T|f(t,\cdot,Y_t,{A^{1/2}}(X_{t-s})Z_t)|dt<\infty, P^\nu-a.s.,$
and
$E^\nu(\int_s^T|{A^{1/2}}(X_{t-s})Z_t|_H^2dt)^{p/2}<\infty,$
and the following equation under $P^\nu$ holds
\begin{equation}Y_t=\xi+\int_t^Tf(r,Y_r,{A^{1/2}}(X_{r-s})Z_r)dr-\int_{t-s}^{T-s}Z_{s+l}.dM_l,\qquad s\leq t\leq T.\end{equation}
The next result gives the probabilistic interpretation of Theorem 3.8. Let us assume that $f:[0,T]\times E\times \mathbb{R}^l\times H^l\rightarrow \mathbb{R}^l$
is the measurable function appearing in the basic equation (3.1).

Let $\phi:E\rightarrow \mathbb{R}^l$ be measurable and for each $p>1$, let $A_p$ denote the set of all points $(s,x)\in[0,T)\times \mathcal{N}^c$ such that
$$E^x\int_s^Tf^{',r}(t,X_{t-s})dt<\infty, \qquad \forall r\geq 0,$$
and
$$E^x(|\phi|^p(X_{T-s})+(\int_s^T|f^0(t,X_{t-s})|ds)^p)<\infty.$$
Set $D:=\cup_{p>1}A_p, A_{p,s}:=\{x\in\mathcal{N}^c,(s,x)\in A_p\}$, and $A_s:=\cup_{p>1}A_{p,s},s\in [0,T)$. By the
same arguments as in [4, Theorem 5.4], we have the following results.
\vskip.10in
\th{Theorem 5.4} Assume that (A4) or (A4') and (A5) holds and  that the function $f$ satisfies conditions (H1),(H2),(H3). Then there exist nearly Borel measurable functions $(u,\psi), u:D\rightarrow \mathbb{R}^l,\psi:D\rightarrow H^l$, such that, for each $s\in[0,T)$ and each $x\in A_{p,s}$, the pair $(u(t,X_{t-s}),\psi(t,X_{t-s}))_{s\leq t\leq T}$ solves  BSDE (5.3) in $L^p(P^x)$ with data $(\phi(X_{T-s}),f(t,X_{t-s},y,z))$ over the interval $[s,T]$.

In particular, the functions $u,\psi$ satisfy the following estimate, for $(s,x)\in A_p$,
$$E^x(\sup_{t\in[s,T]}|u(t,X_{t-s})|^p+(\int_s^T|{A^{1/2}}\psi(t,X_{t-s})|^2dt)^{p/2})\leq KE^x(|\phi(X_{T-s})|^p+(\int_s^T|f^0(t,X_{t-s})|dt)^p).$$
Moreover, if (A1)-(A3) hold and $f$, $\phi$ satisfy  conditions (H4) and (H5), then the complement of $A_{2.s}$ is $\mu$-negligible (i.e. $\mu(A_{2,s}^c)=0$) for each $s\in[0,T)$, the class of $u1_{A_2}$ is an element of $\hat{F}^l$ which is a generalized
 solution of (3.1), $\psi$ represents a version of $\nabla u$ and the following relation holds for each $(s,x)\in D$ and $1\leq i\leq l,$
\begin{equation}u^i(s,x)=E^x(\phi^i(X_{T-s}))+\int_s^TE^xf^i(t,X_{t-s},u(t,X_{t-s}),{A^{1/2}}(X_{t-s})\psi(t,X_{t-s}))dt.\end{equation}

\th{Remark 5.5} As an application we can consider a control problem in the infinite dimensional case by directly extending results in [3, Section 7]. By the results in Section 3 and Section 5, we directly obtain a  mild solution of the Hamilton-Jacobi-Bellman equation. For details see [40].

\section{Examples}

In this section, we give some examples satisfying our assumption (A1)-(A5).
\vskip.10in
\th{Example 6.1} (Ornstein-Uhlenbeck semigroup) Given two separable Hilbert spaces $H$ and $U$, consider the stochastic differential equation
\begin{equation}dX(t)=A_1X(t)dt+BdW(t),\qquad X(0)=x\in H,\end{equation}
where $A_1:D(A_1)\subset H\rightarrow H$ is the infinitesimal generator of a strongly continuous semigroup $e^{tA_1}$, $B: U\rightarrow H$ is a bounded linear operator, and $W$ is a cylindrical Wiener process in $U$. Assume

\no(i) $\|e^{tA_1}\|\leq Me^{\omega t}$ for $\omega<0$, $M\geq 0$, and all $t\geq0$.

\no(ii) For any $t>0$ the linear operator $Q_t$, defined as
$$Q_tx=\int_0^te^{sA_1}Ce^{sA_1^*}xds, \textrm{ } x\in H, t\geq0,$$
where $C=BB^*$,  is of trace class.

\no(iii) $Ce^{tA_1^*}=e^{tA_1}C$.

$\mu$ will denote the Gaussian measure in $H$ with mean $0$ and covariance operator $Q_\infty$. Then the bilinear form $$\mathcal{E}(u,u):=\frac{1}{2}\int_H|C^{1/2}\nabla u|^2d\mu, \textrm{ } u\in\mathcal{F}C_b^\infty,$$
is closable. The closure of $\mathcal{F}C_b^\infty$ with respect to $\mathcal{E}_1$ is denoted by $F$. $(\mathcal{E},F)$ is a generalized Dirichlet form in the sense of Remark 2.1 (iii) with $(E_1,\mathcal{B}(E_1),m)=(H,\mathcal{B}(H),\mu)$, $(\mathcal{A},\mathcal{V})=(\mathcal{E},F)$ and $\Lambda=0$. In particular, it is a symmetric Dirichlet form associated with the O-U process given by (6.1) and
satisfies  conditions (A1)-(A5) (see [13, ChapterII]).

\vskip.10in

 \th{Example 6.2} Let $H$ be a real separable Hilbert space (with scalar product $\langle\cdot,\cdot\rangle$ and norm denoted by $|\cdot|$) and $\mu$  a finite positive measure on $H$. We denote its Borel $\sigma$-algebra by $\mathcal{B}(H)$.
 For $\rho\in L^1_+(H,\mu)$ we consider the following bilinear form
$$\mathcal{E}^\rho(u,v)=\frac{1}{2}\int_H\langle \nabla u,\nabla v\rangle\rho(z)\mu(dz),u,v\in \mathcal{F}C_b^\infty,$$
where $ L^1_+(H,\mu)$ denotes the set of all non-negative elements in $ L^1(H,\mu)$.
 There are many examples for $\rho$ such that $\mathcal{E}^\rho$ is closable. For example if $\rho d\mu$ is a "Log-Concave" measure in the sense of \cite{ASZ}, and more examples can be found in \cite{MR}. The closure of $\mathcal{F}C_b^\infty$ with respect to $\mathcal{E}_1$ is denoted by $F$. $(\mathcal{E},F)$ is a generalized Dirichlet form in the sense of Remark 2.1 (iii) with $(E_1,\mathcal{B}(E_1),m)=(H,\mathcal{B}(H),\mu)$, $(\mathcal{A},\mathcal{V})=(\mathcal{E},F)$ and $\Lambda=0$. In particular, it is a symmetric Dirichlet form and satisfies (A1)-(A5).
  Assume that:
$A_1:D(A_1)\subset H\rightarrow H$ is a linear self-adjoint operator on H such that $\langle A_1x,x\rangle\geq\delta|x|^2 \textrm{ }\forall x\in D(A_1)$ for some $\delta>0$ and $A^{-1}_1$ is of trace class.
 $\mu$ will denote the Gaussian measure in $H$ with mean $0$ and covariance operator$$Q:=\frac{1}{2}A^{-1}_1.$$
 We are concerned with the following two cases.

1. Choose $\rho=\frac{e^{-2U(x)}}{\int_He^{-2U(y)}dy}$ for a Borel map $U:H\rightarrow(-\infty,+\infty]$ with $\int_He^{-2U(y)}dy\in(0,\infty)$. Under some regular condition for $U$,  the process associated with $\mathcal{E}^\rho$ is the solution of the following SPDE
$$dX(t)=(A_1X(t)+\nabla U(X(t))dt+dW(t),\qquad X(0)=x\in H.$$

2. $\rho=1_{\{|x|_H\leq 1\}}$. This case  has been studied in \cite{ASZ}, \cite{RZZ} and it is associated with a reflected O-U process (\cite{RZZ}). The Kolomogorov equation associated with $\mathcal{E}$ has been studied in \cite{BDT} and the solution corresponds to the Kolomogorov equation with Neumman boundary condition.

\vskip.10in
\th{Example 6.3} Consider Example 6.1 and assume that, in addition we are given  a nonlinear  function $F:H\rightarrow H$ such that there exists $K>0$, $|F(x)-F(y)|_H\leq K|x-y|, x,y\in H$ and $\langle F(x)-F(y),x-y\rangle\leq 0, x,y\in H$.  $A_1$ is an operator which satisfies the condition in Example 6.2 and $A^{-1+\delta}_1$ is trace-class for some $\delta\in(0,\frac{1}{2})$. We are concerned with the stochastic differential equation
\begin{equation}dX(t)=(A_1X(t)+F(X(t))dt+BdW(t),\qquad X(0)=x\in H.\end{equation}
The Kolomogrov operator associated with (6.2) is given by
$$K_0\varphi=\frac{1}{2}Tr[CD^2\varphi]+\langle x, A_1^*D\varphi\rangle+\langle F(x),D\varphi\rangle,\varphi\in \mathcal{E}_{A_1}(H),$$
where $\mathcal{E}_{A_1}(H):=\textrm{ linear span } \{\varphi_h(x)=e^{i\langle h,x\rangle}: h\in D(A^*_1)\}$.
Assume the semigroup $e^{tA_1}$ is analytic. Then by [15, Theorem 11.2.21] there exists a unique invariant measure $\nu$ for $K_0$ i.e.
$$\int K_0\varphi d\nu=0, \textrm{ for all }\varphi\in\mathcal{E}_{A_1}(H),$$
and $\nu$ is absolutely continuous with respect to $\mu$ from Example 6.1 and for $\rho=\frac{d\nu}{d\mu}$ we have that $\rho\in W^{1,2}(H,\mu)$ and $D\log\rho\in W^{1,2}(H,\nu;H).$

 Then by  [38, Section 4.2], we know that the bilinear form on $L^2(H;\nu)$ associated with $K_0$ is a generalized Dirichlet form in the sense of Remark 2.1 (iii) with $(E_1,\mathcal{B}(E_1),m)=(H,\mathcal{B}(H),\nu)$, $(\mathcal{A},\mathcal{V})=(0,L^2(H,\nu))$ and $\Lambda=K_0$. It
  satisfies conditions (A1)-(A5). There are even more general conditions on $F$ and $A_1$ which can be found in [15, Theorem 5.2] such that conditions (A1)-(A5) hold.
\vskip.10in

The following example is given in [38, Section 4.2].

\th{Example 6.4} Assume that $E$ is a separable real Hilbert space with inner product $\|\cdot\|_E^{1/2}$ and $H\subset E$ densely by a Hilbert-Schmidt map. Let $B:E\rightarrow E$ be a Borel measurable vector field satisfying the following conditions:

(B.1) $\lim_{\|z\|_E\rightarrow\infty}\langle B(z),z\rangle=-\infty.$

(B.2) $_{E'}\langle l, B\rangle_E:E\rightarrow\mathbb{R}$ is weakly continuous for all $l\in E'$.

(B.3) There exist $C_1,C_2,d\in (0,\infty)$, such that $\|B(z)\|_E\leq C_1+C_2\|z\|_E^d$.

Then by [9, Theorem 5.2] there exists a probability measure $\mu$ on $(E,\mathcal{B}(E))$ such that $_{E'}\langle l, B\rangle_E\in L^2(E;\mu)$ for all $l\in E'$ and such that
$$\int\frac{1}{2}\Delta_H u+\frac{1}{2}   _{E'}\langle\nabla u,B\rangle_Ed\mu=0 \textrm{ for all } u\in \mathcal{F}C_b^\infty,$$
where $\Delta_H$ is the Gross-Laplacian, i.e., $\Delta_Hu=\sum_{i,j=1}^m\frac{\partial f}{\partial x_i\partial x_j}(l_1(z),.., l_m(z))\langle l_i,l_j\rangle_H$ for $u=f(l_1,...,l_m)\in \mathcal{F}C_b^\infty$. Assume $B(z)=-z+v(z), v:E\rightarrow H$.
For the bilinear form associated with  $Lu=\frac{1}{2}\Delta_H u+\frac{1}{2}_{E'}\langle \nabla u, B\rangle_E, u\in \mathcal{F}C_b^\infty$ on $L^2(E,\mu)$ is a generalized Dirichlet form in the sense of Remark 2.1 (iii) with $(E_1,\mathcal{B}(E_1),m)=(H,\mathcal{B}(H),\nu)$, $(\mathcal{A},\mathcal{V})=(0,L^2(H,\nu))$ and $\Lambda=L$. It satisfies conditions (A1)-(A5).

\section{Appendix}

\th{Appendix A. Basic Relations for the Linear Equation}
In this section we assume that (A1)-(A3) hold.

\vskip.10in
\th{Lemma A.1} If $u$ is a bounded generalized solution of equation (2.3), then $u^+$ satisfies the following relation with $0\leq t_1<t_2\leq T$
$$\|u_{t_1}^+\|_2^2\leq 2\int_{t_1}^{t_2}(f_s,u_s^+)ds+\|u_{t_2}^+\|_2^2.$$

\proof Choose  the approximation sequence $u^n$ for $u$ as in the existence proof of Proposition 2.7. Denote its related data as $f^n,\phi^n$ .

Suppose that  the following holds
$$\|(u_{t_1}^n)^+\|_2^2\leq 2\int_{t_1}^{t_2}(f_s^n,(u_s^n)^+)ds+\|(u_{t_2}^n)^+\|_2^2,\eqno(A.1)$$
where $0\leq t_1\leq t_2\leq T$. Since
$\|u^n\|_2$ are uniformly bounded, we have $\lim_{n\rightarrow\infty}\int_{t_1}^{t_2}(f_s^n,(u_s^n)^+)ds=\int_{t_1}^{t_2}(f_s,u_s^+)ds.$
By letting $n\rightarrow\infty$  in equation (A.1) the assertion follows.

Therefore, the problem is reduced to the case where $u$ belongs to $b\mathcal{C}_T$ ; in the
remainder we assume $u\in b\mathcal{C}_T$.  (2.7), written with $u^+\in bW^{1,2}([0,T];L^2)\cap L^2([0,T];F)$ as test functions, takes the form
$$\aligned& \int_{t_1}^{t_2}(u_t,\partial_t(u_t^+))dt+\int_{t_1}^{t_2}\mathcal{E}^{A}(u_t,u_t^+)dt+\int_{t_1}^{t_2}\int\langle A^{1/2}b, D_{A^{1/2}} u_t
\rangle u_t^+dmdt
\\=&\int_{t_1}^{t_2}(f_t,u_t^+)dt+(u_{t_2},u^+_{t_2})-(u_{t_1},u^+_{t_1}).\endaligned\eqno(A.2)$$
By [5, Theorem 1.19], we  obtain
$\int_{t_1}^{t_2}(u_t,\partial_t(u_t^+))dt=\frac{1}{2}(\|u_{t_2}^+\|_2^2-\|u_{t_1}^+\|_2^2).$
Then
$$\aligned&\|u_{t_1}^+\|_2^2 +2\int_{t_1}^{t_2}\mathcal{E}^A(u_t,u_t^+)dt+2\int_{t_1}^{t_2}\int\langle A^{1/2}b, D_{A^{1/2}} u_t
\rangle_H u_t^+dmdt
=2\int_{t_1}^{t_2}(f_t,u_t^+)dt+\|u_{t_2}^+\|_2^2.\endaligned\eqno(A.3)$$
Next we prove for $u\in bF$
$$\mathcal{E}(u,u^+)\geq0.\eqno(A.4)$$
We have the above relation for $u\in \mathcal{D}(L)$. For $u\in bF$, by (A3) we  choose a uniformly bounded sequence $\{u_n\}\subset \mathcal{D}(L)$ such that $\mathcal{E}^{A}_{1}(u_n-u)\rightarrow0$.  And then we have
$$|\int\langle A^{1/2}b, D_{A^{1/2}} u
\rangle_H u^+d\mu-\int\langle A^{1/2}b, D_{A^{1/2}} u_n
\rangle_H u_n^+d\mu|
\rightarrow0. $$
Since $\mathcal{E}^A(u^+)\leq \mathcal{E}^A(u)$, $\sup_n\mathcal{E}^A(u^+_n)\leq \sup_n\mathcal{E}^A(u_n)<\infty$, we also have
$|\mathcal{E}^A(u_n,(u_n)^+)-\mathcal{E}^A(u,u^+)|
\rightarrow0.$
As a result we have (A.4) for bounded $u\in F$.
So we have $\|u_{t_1}^+\|_2^2 \leq2\int_{t_1}^{t_2}(f_t,u_t^+)dt+\|u_{t_2}^+\|_2^2.$
$\hfill\Box$

\vskip.10in
To extend the class of solutions we are working with,  to allow $f$ to belong to $L^1(dt\times d\mu)$, we need the following proposition. It is a modified version of the above lemma and proof is similar as before.

\vskip.10in
\th{Lemma A.2} Let $u\in\hat{F}$ be bounded and $f\in L^1(dt\times d\mu)$, be such that the weak relation (2.7) is satisfied with  test functions in $b\mathcal{C}_T$ and some function $\phi\geq0$, $\phi\in L^2\cap L^\infty$. Then $u^+$ satisfies the following relation for
$0\leq t_1<t_2\leq T$
$$\|u_{t_1}^+\|_2^2\leq 2\int_{t_1}^{t_2}(f_s,u_s^+)ds+\|u_{t_2}^+\|_2^2.$$

\vskip.10in

The next proposition is a modification of [4, Proposition 2.9]. It represents a version of the maximum principle.

\vskip.10in
\th{Proposition A.3} Let $u\in\hat{F}$ be bounded and $f\in L^1(dt\times d\mu),f\geq0$, be such that the weak relation (2.7) is satisfied with test functions in $b\mathcal{C}_T$ and some function $\phi\geq0$, $\phi\in L^2\cap L^\infty$. Then $u\geq 0$ and it is represented by the following relation:
$$ u_t=P_{T-t}\phi+\int_t^TP_{s-t}f_sds.$$
Here we use $P_t$ is a $C_0$-semigroup on $L^1(E;\mu)$ to make $P_{s-t}f_s$ meaningful.

\proof Let $(f^n)_{n\in N}$ be a sequence of bounded functions on $[0,T]\times E$ such that
$0\leq f^n\leq f^{n+1}\leq f, \lim_{n\rightarrow\infty}f^n=f.$
Since $f^n$ is bounded, we have $f^n\in L^1([0,T];L^2)$. Next we define
$ u_t^n=P_{T-t}\phi+\int_t^TP_{s-t}f_s^nds.$
Then $u^n\in\hat{F}$ is the unique generalized solution for the data $(\phi,f^n)$. Clearly $0\leq u^n\leq u^{n+1}$ for $n\in N$. Define $y:=u^n-u$ and $\tilde{f}:=f^n-f$. Then $\tilde{f}\leq 0$ and $y$ satisfies the weak relation (2.7) for the data $(0,\tilde{f})$. Therefore by Lemma A.2, we have for $t_1\in [0,T]$,
$\|y_{t_1}^+\|_2^2\leq 2\int_{t_1}^{T}(\tilde{f}_s,y_s^+)ds\leq 0.$
We conclude that $\|y_{t_1}^+\|_2^2=0$. Therefore, $u\geq u^n\geq0$ for $n\in N$. Set $v:=\lim_{n\rightarrow\infty}u^n$. And we  obtain that
 $\lim_{n\rightarrow\infty}\|u_t^n-v_t\|_2^2=0$ and
$\lim_{n\rightarrow\infty}|\int_t^T\int(f_s^nu_s^n-f_sv_s)d\mu ds|=0.$
By [25, Lemma 2.12] we have
$$\int_t^T\mathcal{E}^A(v_s)ds\leq \int_t^T\liminf_{n\rightarrow\infty}\mathcal{E}^A(u_s^n)ds\leq\liminf_{n\rightarrow\infty}\int_t^T\mathcal{E}^A(u_s^n)ds.$$
Finally, we get for $t\in[0,T]$
$$\aligned&\|v_t\|_2^2+2\int_t^T\mathcal{E}^A(v_s)ds\leq\lim_{n\rightarrow\infty}\|u_t^n\|_2^2+
2\liminf_{n\rightarrow\infty}\int_t^T\mathcal{E}^A(u_s^n)ds\\
\leq &\lim_{n\rightarrow\infty}(2\int_t^T(f_s^n,u_s^n)ds+\|\phi\|_2^2+2\alpha\int_t^T\|u_s^n\|_2ds)
=2\int_t^T(f_s,v_s)ds+\|\phi\|_2^2+2\alpha\int_t^T\|v_s\|_2ds.\endaligned$$
Since the right side of this inequality is finite and $t\mapsto v_t$ is $L^2$-continuous, it follows that $v\in\hat{F}$.

Now we show that $v$ satisfies the weak relation (2.7) for the data $(\phi,f)$. As $\varphi^n(t):=\|u_t^n-v_t\|_2$ is continuous and decreasing to $0$, we conclude by Dini's theorem
$\lim_{n\rightarrow\infty}\sup_{t\in[0,T]}\|u_t^n-v_t\|_2=0,$
and therefore
$\lim_{n\rightarrow\infty}\int_0^T\|u_t^n-v_t\|_2^2=0.$
Furthermore, there exists $K\in \mathbb{R}_+$ and a subsequence $(n_k)_{k\in N}$ such that
$|\int_0^T\mathcal{E}^A(u_s^{n_k})ds|\leq K , \forall k\in N.$
In particular, $\int_0^T\int|D_{A^{1/2}} u_s^{n_k}|^2_Hd\mu ds\leq  K , \forall k\in N.$
We obtain $\lim_{k\rightarrow\infty}\int_0^T\mathcal{E}^A(u_s^{n_k},\varphi_s)ds=\int_0^T\mathcal{E}^A(v_s,\varphi_s)ds,$
and
$$\lim_{k\rightarrow\infty}\int_0^T\int\langle A^{1/2}b, D_{A^{1/2}} u_s^{n_k}\rangle_H\varphi_sd\mu ds=\int_0^T\int\langle A^{1/2}b, D_{A^{1/2}} v_s\rangle_H\varphi_sd\mu ds,$$
 which implies (2.7) for $v$ associated to $(\phi,f)$. Clearly $u-v$ satisfies (2.7) with data $(0,0)$ for $\varphi\in b\mathcal{C}_T$. By Proposition 2.7 we have $u-v=0$. Since
$ v_t=P_{T-t}\phi+\int_t^TP_{s-t}f_sds,$
the assertion follows.$\hfill\Box$

\vskip.10in
\th{Corollary A.4} Let $u\in\hat{F}$ be bounded and $f\in L^1(dt\times d\mu)$ be such that the weak relation (2.7) is satisfied with test functions in $b\mathcal{C}_T$ and some function  $\phi\in L^2\cap L^\infty$. Assume there exists $g\in bL^1(dt\times d\mu)$ such that $f\leq g$. Then $u$  has the following representation:
$$ u_t=P_{T-t}\phi+\int_t^TP_{s-t}f_sds.$$
\proof Define $f^n:=(f\vee(-n))\wedge g, n\in \mathbb{N}$. Then $(f^n)_{n\in \mathbb{N}}$ is a sequence of bounded functions such that $f^n\downarrow f$ and $f^n\leq g$ then by the same arguments as in Proposition A.3, the assertion follows.$\hfill\Box$

By the above results and Proposition 2.8, we obtain the following lemma by the same arguments as the proof
in [4, Lemma 2.12].
\vskip.10in

\th{Lemma A.5} If $f,g\in L^1([0,T];L^2)$ and $\phi\in L^2$, then the following relations hold $\mu$-a.e.:
$$\int_t^TP_{s-t}(f_sP_{T-s}\phi)ds\leq \frac{1}{2}P_{T-t}\phi^2+\int_t^T\int_s^TP_{s-t}(f_sP_{r-s}f_r)drds.\eqno(A.5)$$

\vskip.10in
\th{Acknowledgement.} I would like to thank Professor M. R\"{o}ckner for valuable discussions and for suggesting me to associate BSDE's to generalized Dirichlet forms which was one motivation for this paper . I would also like to thank Professor Ma Zhiming and Zhu Xiangchan for their helpful discussions.

\end{document}